%% file: syz.tex
\documentclass[12pt]{amsart}
\usepackage{amscd,amsthm,amssymb,amsfonts,amsmath,epsfig,epic,eclbip}
\usepackage{bbm}
\usepackage{url}
\usepackage{color}

\theoremstyle{plain}
\newtheorem{thm}{Theorem}[section]
\newtheorem{lemma}[thm]{Lemma}
\newtheorem{prop}[thm]{Proposition}

\theoremstyle{definition}
\newtheorem*{defn}{Definition}

\theoremstyle{remark}
\newtheorem*{remark}{Remark}
\newtheorem*{ack}{Acknowledgments}

\newcommand{\Z}{\mathbb Z}    
\newcommand{\R}{\mathbb R}    
\newcommand{\C}{\mathbb C}    
\newcommand{\PP}{\mathbb P}   
\newcommand{\SL}{\operatorname{SL}}

\newcommand{\<}{\langle}   
\renewcommand{\>}{\rangle} 

\newcommand{\D}{\Delta}
\newcommand{\DSv}{\Delta^\vee_\lambda}
\newcommand{\Dv}{\Delta^\vee}
\newcommand{\DT}{{\Delta_\nu}}
\newcommand{\smooth}{\Sigma\backslash D}
\newcommand{\carrier}{\operatorname{carrier}}
\newcommand{\bsd}{\operatorname{bsd}}
\newcommand{\suchthat}{\ : \ }

\newcommand{\dD}{{\partial\Delta}}
\newcommand{\dDv}{{\partial\Delta^\vee}}
\newcommand{\U}{\mathcal U}
\newcommand{\V}{\mathcal V}
\newcommand{\T}{\mathbb T}
\newcommand{\base}{\Sigma\backslash N(D)}
\newcommand{\Log}{\mathrm{Log}}
\newcommand{\am}{\mathcal A^\lambda}

\newcommand{\aff}{{\operatorname{af{}f}}}
\newcommand{\join}{\ast}
\renewcommand{\setminus}{\backslash}
\newcommand{\Star}{\operatorname{star}}
\newcommand{\vertS}{\operatorname{vert}(S)}
\newcommand{\vertT}{\operatorname{vert}(T)}

\newcommand{\relint}{\operatorname{relint}}
\newcommand{\conv}{\operatorname{conv}}
\newcommand{\pull}{\operatorname{pull}}
\newcommand{\dDT}{\partial \DT}
\renewcommand{\L}{\mathcal L}
\newcommand{\chain}{\bsd}

\newcommand{\moment}{\mu}
\newcommand{\vol}{\operatorname{vol}}
\newcommand{\Hsm}{H_s^{\mathrm{sm}}}

\newcommand{\dDS}{\partial \DSv} 
\newcommand{\dU}{\partial\U}
\newcommand{\dV}{\partial\V}

\newcommand{\Ve}{V^\delta}
\newcommand{\Ue}{U^\epsilon}
\newcommand{\We}{W^\epsilon}
\newcommand{\Qev}{Q^\lambda_{(\{0\}\mid v)}(\epsilon)}
\newcommand{\Qew}{Q^\lambda_{(\{0\}\mid w^\perp)}(\epsilon)}
\newcommand{\Qv}{Q^\lambda_{(\{0\}\mid v)}(0)}
\newcommand{\Qw}{Q^\lambda_{(\{0\}\mid w^\perp)}(0)}
\newcommand{\QeIJ}{Q^\lambda_{(I\mid J)}(\epsilon)}
\newcommand{\QIJ}{Q^\lambda_{(I\mid J)}(0)}
\newcommand{\ball}{Q^\lambda_{\{0\}}(\epsilon)}
\newcommand{\Cstar}{\C\setminus\{0\}}
\newcommand{\cone}{\operatorname{cone}}

\newcommand{\F}{\mathcal F}
\newcommand{\X}{\mathfrak X}

\renewcommand{\caption}[1]{%
  \refstepcounter{figure}
  \begin{minipage}[t]{.8\textwidth}
    \vspace{.3\baselineskip}
    \begin{center}
      \scriptsize {\sc Figure \thefigure:}
      #1
    \end{center}
  \end{minipage}
}

\begin{document}

\begin{abstract}
  We describe in purely combinatorial terms dual pairs of integral
  affine structures on spheres which come from the conjectural metric
  collapse of mirror families of Calabi-Yau toric hypersurfaces. The
  same structures arise on the base of a special Lagrangian torus
  fibration in the Strominger-Yau-Zaslow conjecture. We study the
  topological torus fibration in the large complex structure limit
  and show that it coincides with our combinatorial model. 
\end{abstract}

\title[Integral affine structures on spheres I]{Integral affine
  structures on spheres \makebox[0pt]{\raisebox{1in}{\normalfont
      DUKE-CGTP-02-05}} \\ and torus fibrations of Calabi-Yau toric
  hypersurfaces I } 
\author{Christian
Haase} 
\author{Ilia Zharkov}
\address{Mathematics Department \\ Duke University \\ Durham, NC 27708
  \\ USA}
\email{[haase,zharkov]@math.duke.edu}
\maketitle
\section{Introduction} 
\input{intro}
\begin{ack}
  We are indebted to David Morrison for the original idea and
  continuous suggestions throughout our work.
  We thank Anda Degeratu for valuable conversations, and the Duke
  Math/Physics group for a stimulating environment. The second author  
  would also like to thank IHES where he stayed during the final stage of 
  the work, for its hospitality and financial support.
\end{ack}
\section{The combinatorial model} \label{sec:model}
\input{model}
\subsection{$\SL(n,\Z)$-structure and monodromy}
\input{monodromy}
\section{Torus fibrations of Calabi-Yau toric hypersurfaces}
\label{sec:hypersurfaces}
\subsection{The family}
\input{family}
\subsection{Amoebas of hypersurfaces}
\input{moment1}
\subsection{Foliation of $\R^d\setminus\ball$}
\input{moment2}
\subsection{The torus fibration}
\input{moment4}
\section{Outlook}
\input{outlook1}
\input{outlook2}
\bibliographystyle{alpha}
\bibliography{references}

\end{document}

%% file: intro.tex
Since Strominger, Yau, and Zaslow conjectured an interpretation of
mirror symmetry as a duality of special Lagrangian torus
fibrations \cite{SYZ}, there has been considerable progress toward
proving the topological consequences of this conjecture. Recently,
for a Calabi-Yau family in a neighborhood of the large complex
structure point, Kontsevich and Soibelman \cite{KS} and Gross and
Wilson \cite{GW} conjectured the existence of an affine K\"ahler
structure on the (complement of a codimension two locus of the)
limiting metric space of the Gromov-Hausdorff collapse. This metric
space should be identified with the base of the SYZ fibration (with
McLean's metric). In particular, the collapse picture asserts the
existence of an integral affine structure on the base,
which is enough to reconstruct the non-degenerate part of the
topological torus fibration of the CY family.

An integral affine structure on the base of the conjectural SYZ torus
fibration has been described by Ruan \cite{Ru} and Gross \cite{Gr} in
the quintic case, followed by work of Ruan \cite{Ru2} for general
toric hypersurfaces. The present paper provides a short and
explicit combinatorial description of an affine structure with mirror
duality built in, based on an idea of David
Morrison~\cite{MorrisonLecture}.

Given a dual pair of $d$-dimensional reflexive polytopes 
with coherent triangulations of their boundaries,
we construct in \S\ref{sec:SigmaD} a $(d-1)$-dimensional polytopal
complex $\Sigma$, topologically a sphere, with a codimension 2
subcomplex $D$, which we call a discriminant locus. The manifold
$\Sigma \backslash D$ possesses an integral affine structure
(\S\ref{sec:monodromy}). That is, the tangent space at any point in
$\smooth$ contains a natural integral lattice, and one can form a
$(d-1)$-torus fibration $W\to\smooth$ by taking fiber wise
quotients. The nerve of the covering of $\smooth$ by affine charts is
a two-colored graph, whose nodes are labeled by the vertices
of the triangulations and an edge connects any two which live in dual
faces of the polytopes.

The more technical \S\ref{sec:bsd} is
concerned with the sphericity of $\Sigma$. We develop a generalization
of barycentric subdivisions which might be of independent interest to
a combinatorially inclined reader.

In Section~\ref{sec:hypersurfaces} we link the model to the topology of 
toric hypersurfaces $H_s$ constructed from our input data
(\S\ref{sec:family}). The main result Theorem~\ref{thm:main} asserts
that for any neighborhood $N$ of $D$, and a hypersurface with large
enough complex structure, there is a torus fibration of $\Hsm$, a
portion of the hypersurface, over $\Sigma \setminus N$, which is
diffeomorphic to the restriction of our model fibration.
%
%

In the second part of the paper we will develop a connection of
our model to the {\it geometry} of the hypersurfaces, conjectured in
\cite{GW} and \cite{KS}. Namely, we will show that the above
diffeomorphism provides, in fact, an ``almost'' holomorphic embedding
(there is a preferred choice of complex structures on $W$). Moreover,
we will construct a family of K\"ahler forms on $H_s$ in the expected
class, so that the pairs $(H_s,H_s\setminus\Hsm)$ with the induced
metrics converge in the Gromov-Hausdorff sense to the pair
$(\Sigma,D)$.


%% file: model.tex
We are given a reflexive polytope $\Dv \subset \R^d$. That is, $\Dv$
is the convex hull of finitely many points in the lattice $\Z^d$, it
contains the origin in the interior, and the vertices of the dual
polytope $\D = \{ m \in (\R^d)^* \suchthat \<m,n\> \le 1$ for all $n
\in \Dv \}$ belong to the dual lattice
$(\Z^d)^*$~\cite{BatyrevMirror}.

Another part of our input are two sufficiently generic vectors
$\lambda \in \Z^{\D \cap (\Z^d)^*}$, and $\nu \in \Z^{\Dv \cap
  \Z^d}$, which induce central coherent triangulations of $\D$ and
$\Dv$. (Cf.~\cite{LeeHandBook} for coherent/regular triangulations.)
These triangulations restrict to triangulations $S$ and $T$ on the
boundaries $\dD$ and $\dDv$.
\vspace{\baselineskip}
\begin{figure}[htbp]
  \begin{center}
    \input{delta.pstex_t}        \label{fig:DandDv}
    \qquad \qquad
    \input{deltaV.pstex_t}
    \caption{$\D = \conv \left[ 
        \begin{smallmatrix}
          0&\ 2&-2&\ 0&\ 0 \\
          0&\ 0&\ 0&\ 2&-2 \\
          1&-1&-1&-1&-1
        \end{smallmatrix} \right]$
        \quad
        $\Dv = \conv \left[ 
          \begin{smallmatrix}
            \ 0&1&\ 1&-1&-1\\
            \ 0&1&-1&\ 1&-1\\
            -1&1&\ 1&\ 1&\ 1
          \end{smallmatrix} \right]$
        \\[.5\baselineskip]
        The values of $\lambda$ and $\nu$ are marked on the vertices,
        $\lambda(0) = \nu(0) = 0$.}
  \end{center}
\end{figure}
\vspace{\baselineskip}
Using $\lambda$ and $\nu$ we define polytopes
\begin{align*}
  \DT &= \{ m \in (\R^d)^* \suchthat \<m,n\> \le \nu(0)-\nu(n) \text{
    for all } n \in \Dv \cap \Z^d \} \\
  \DSv &= \{ n \in \R^d \suchthat \<m,n\> \le \lambda(0)-\lambda(m)
  \text{ for all } m \in \D \cap (\Z^d)^* \} ,
\end{align*}
whose normal fans are given by the cones spanned by the faces of $T$
respectively $S$.
\begin{figure}[htbp]
  \begin{center}
    \includegraphics{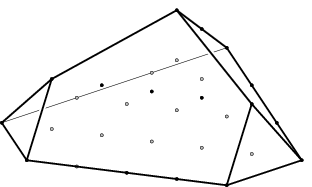}
    \qquad \qquad
    \includegraphics{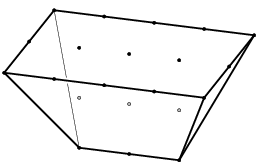}
    \caption{\quad $\DT$ \qquad \qquad \qquad \qquad \qquad \qquad
      \qquad $\DSv$ \qquad \qquad \qquad \qquad }
  \end{center}
\end{figure}
\subsection{The base and the discriminant locus}
\label{sec:SigmaD}
Consider the product $\D \times \Dv \subset (\R^d)^* \times \R^d$. The
complex $\Sigma$ --- our prospective base space --- will be a
subdivision of
\begin{equation} \label{eq:sphere}
  |\Sigma| = \{ (m,n) \in \D \times \Dv \suchthat \< m,n \> = 1 \}
\end{equation}
which is the union of all sets of the form $F\times F^\vee$, where $F$
runs over all proper faces of $\D$, and $F^\vee$ denotes the face (of
$\Dv$) dual to $F$.
The triangulations $S$ and $T$ induce a subdivision of $\D \times
\Dv$ into products of simplices, which restricts to a subdivision of
$|\Sigma|$.
\begin{figure}[htbp]
  \begin{center}
    \includegraphics{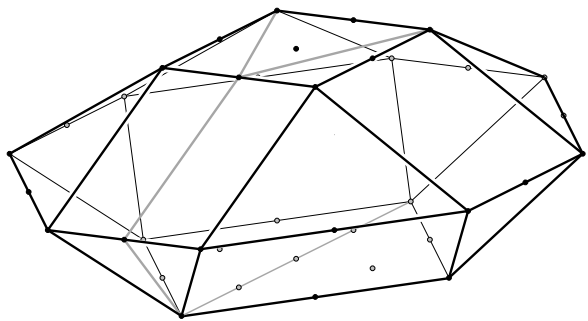}
    \caption{ The subdivision $S \times T$ of $|\Sigma|$ into products
      of simplices.}
  \end{center}
\end{figure}
\begin{defn}
  $\Sigma$ is the restriction to $|\Sigma|$ of the product
  subdivision $\bsd(S) \times \bsd(T)$ of $\D \times \Dv$.
\end{defn}
The vertices of $\Sigma$ correspond to pairs
$(\widehat{\sigma},\widehat{\tau})$ of (barycenters of) simplices $\sigma \in
S$ and $\tau \in T$ with $\< \sigma,\tau \> = 1$. Under this
correspondence, the cells of $\Sigma$ correspond to pairs of chains in
the face posets of $S$ and $T$.
\begin{lemma}\label{lemma:sphere}
  $\Sigma$ is isomorphic to the boundary complex of a $d$-dimensional
  polytope, and therefore topologically a $(d-1)$-sphere.
\end{lemma}
We postpone the proof, and proceed with definitions.
\begin{defn}
  The singular locus $D$ is the full subcomplex of $\Sigma$, induced
  by vertices $(\widehat{\sigma},\widehat{\tau})$, such that neither
  $\sigma$ nor $\tau$ is $0$-dimensional. (This set indeed induces a
  subcomplex.)
\end{defn}
For $d=2$, $D$ is empty. For $d=3$, $D$ will be a finite set of $\le
24$ points (with equality if $S$ and $T$ use all lattice points). For
$d=4$, $D$ will be the first subdivision of a trivalent graph.
%
For all $d$, the complement $\smooth$ is homotopy equivalent to a
bipartite graph.
\begin{defn}
  Let $\Gamma$ be the graph on the vertex set $\vertS \cup \vertT$
  with an edge between $v \in \vertS$ and $w \in \vertT$ if and only
  if $\< v,w \> = 1$.
\end{defn}
\begin{lemma} \label{lemma:complement}
  $\smooth$ is homotopy equivalent to $\Gamma$.
\end{lemma}
First, we introduce some notation that will be used later on.
There are two piecewise linear projections
\begin{equation*}
  p_1 \colon \Sigma \rightarrow \bsd(S) \ \text{ and } \ 
  p_2 \colon \Sigma \rightarrow \bsd(T) .
\end{equation*}
For a vertex $v \in \vertS$ or $w \in \vertT$, define $U_v$
respectively $V_w$ to be the preimages
\begin{equation*}
  U_v = p_1^{-1}(\Star_{\bsd(S)}(v)) \ \text{ and } \ 
  V_w = p_2^{-1}(\Star_{\bsd(T)}(w))
\end{equation*}
of open stars in the barycentric subdivisions. Thus, $U_v$ is an open
regular neighborhood of the contractible set $v \times
(\carrier_\D v)^\vee$ in $\Sigma$, and $V_w$ is an open regular
neighborhood of $(\carrier_{\Dv}w)^\vee \times w$.
We will abbreviate the collections by $\U = (U_v)_{v \in \vertS}$ and
$\V = (V_w)_{w \in \vertT}$. For future reference, define
$\overline{\U} = (\overline{U}_v)_{v \in \vertS}$ and $\overline{\V} =
(\overline{V}_w)_{w \in \vertT}$, as well as $\dU = \bigcup \partial 
U_v$ and $\dV = \bigcup \partial V_w$. 
Observe that with these definitions $\bigcup \U \cup \bigcup \V =
\smooth$, and $\dU \cap \dV = D$.
\begin{figure}[htbp]
  \begin{center}
    \includegraphics{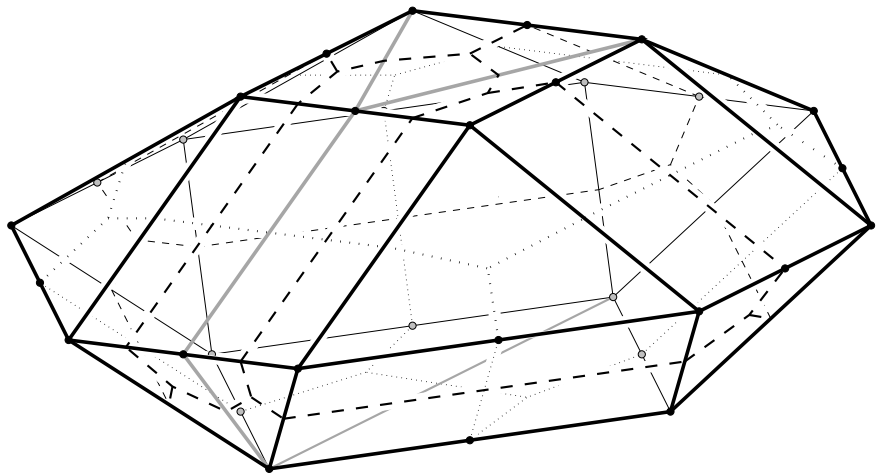}
  \caption{ The dotted lines are $\dU$, and the dashed lines are
  $\dV$. Their intersection $D$ consists of $9$ points.}
\end{center}
\end{figure}
\begin{proof}[Proof of Lemma~\ref{lemma:complement}]
  Two members $U_v$ and $V_w$ of this covering intersect if and only
  if
  \begin{equation*}
    v \in (\carrier_{\Dv} w)^\vee \iff w \in (\carrier_\D v)^\vee \iff
    \< v,w \> = 1 ,
  \end{equation*}
  in which case they intersect in the contractible set
  $\Star_\Sigma((v,w))$. So the claimed homotopy equivalence follows
  from the nerve lemma (cf.~e.g.\ \cite[Thm.~10.6]{Bjoerner}).
\end{proof}
\subsection{The proof of Lemma~\ref{lemma:sphere}} \label{sec:bsd}
We will describe a coherent subdivisions of $\DSv$ (alternatively
$\DT$), which is isomorphic to $\Sigma$. We hope that the method,
which generalizes barycentric subdivisions, may find other
applications in the future.
\begin{defn}
  Suppose that $P \subset Q \subset \R^d$ are polytopes such that $P$
  is contained in the relative interior of some face of $Q$. Then the
  result of pulling $P$ is by definition the coherent subdivision
  $\pull_P(Q)$ of $Q$ which is induced by the heights $1$ on $P$, and
  $0$ on all faces of $Q$ that do not contain $P$.
\end{defn}
The case $P \subset \relint(Q)$ is used in~\cite{GoodmanPach} in order
to prove that the space between $P$ and $Q$ can be triangulated
without new vertices.
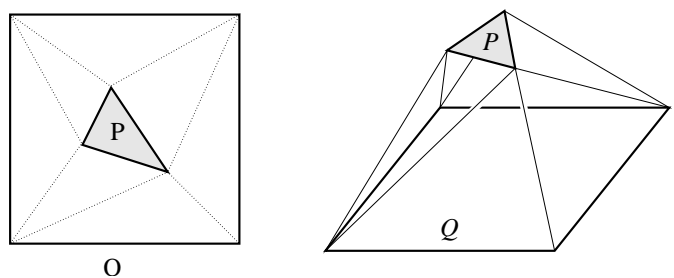
\begin{figure}[htbp]
  \begin{center}
    \input{pull.pstex_t}
    \caption{The pulling subdivision $\pull_P(Q)$.}
  \end{center}
\end{figure}
We can describe this subdivision combinatorially as follows. The
faces of $\pull_P(Q)$ are all sets of the form $\conv(F \cup F')$ for
faces $F \preceq P$ and $F' \prec Q$ with $F \not \subset F'$ (including
$F' = \emptyset$) whose normal cones intersect in their relative
interiors. That is, there should be a normal vector $n \in (\R^d)^*$
which is maximized over $P$ precisely on $F$, and over $Q$ precisely
on $F'$. More generally, if $P$ is in the relative interior of a face
of a polyhedral complex, then pulling $P$ will affect all faces that
contain $P$ in the manner outlined above. If the original subdivision
was coherent, then the pulling subdivision will remain coherent.
This procedure generalizes the well studied pulling subdivisions where
$P$ is a point~\cite[\S~5.2]{Gruenbaum} or \cite{LeeHandBook}.

We use these pullings in order to construct a generalized
barycentric subdivisions below. We start with a purely combinatorial
definition. For a poset $\mathcal{Q}$, the poset/simplicial complex of
chains in $\mathcal{Q}$ is denoted $\chain(\mathcal{Q})$.
\begin{defn}
  Suppose $\kappa \colon \mathcal{Q} \rightarrow \mathcal{P}$ is an
  order preserving, non-rank-increasing correspondence between the
  graded posets $\mathcal{Q}$ and $\mathcal{P}$. Define the
  barycentric subdivision $\bsd(\mathcal{Q},\kappa)$ of $\mathcal{Q}$
  with respect to $\kappa$ as the subposet
  \begin{equation*}
    \{ ( p \ , \ q_0 \prec q_1 \prec \ldots \prec q_r ) \in
    \mathcal{P} \times \chain(\mathcal{Q}) \suchthat  p \preceq
    \kappa(q_0) \}
  \end{equation*}
  of the product poset $\chain(\mathcal{Q}) \times \mathcal{P}$.
\end{defn}
If $\mathcal{P}$ has only one element, this specializes to
$\chain(\mathcal{Q})$. In our applications, $\kappa$ will be clear
from the context, and will write $\bsd(\mathcal{Q},\mathcal{P})$
instead.

One example of such a $\kappa$ arises in the following situation.
Say that a polytope $P \subset \R^d$ is a Minkowski summand of the
polytope $Q \subset \R^d$, if there is an $\varepsilon > 0$, and a
polytope $P' \subset \R^d$ such that $\varepsilon P + P' = Q$. This is
true if and only if the normal fan of $Q$ refines the normal fan of
$P$~\cite{Smilansky}. So there is an order preserving,
non-rank-decreasing correspondence on the level of the normal fans
which turns into $\kappa \colon \L(Q) \rightarrow \L(P)$ on the level
of the face lattices.
\begin{defn} \label{def:bsd}
  Suppose that $P$ is a Minkowski summand of $Q$. Define the
  barycentric subdivision $\bsd(Q;P)$ of $Q$ with respect to $P$ as
  follows.
  Start with $F=Q$, and proceed by decreasing dimension of faces $F
  \preceq Q$. Pull the translate $\varepsilon \kappa(F) - \varepsilon
  \widehat{\kappa(F)} + \widehat{F}$ of the corresponding face of
  $\varepsilon P$ in the relative interior of $F$.
\end{defn}
The usual barycentric subdivision appears as the special case where
$P$ is a point.
\begin{figure}[htbp]
  \begin{center}
    \includegraphics{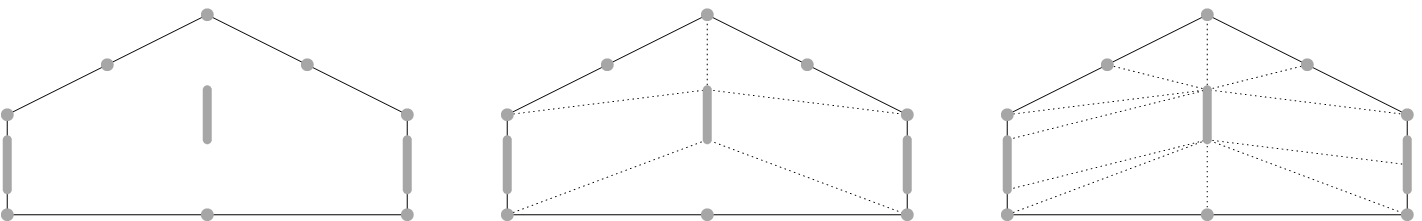}
    \caption{The barycentric subdivision of $\conv \left[ 
        \begin{smallmatrix}
          0&1&\ 1&\ 2&\ 2\\
          0&2&-2&\ 2&-2\\
          1&0&\ 0&-1&-1
        \end{smallmatrix} \right] \prec \DT$
      with respect to the corresponding face $\conv \left[ 
        \begin{smallmatrix}
          0&\ 2\\0&\ 0\\0&-1
        \end{smallmatrix} \right] \prec \D$.}
  \end{center}
\end{figure}
Combinatorially, the face lattice $\L(\bsd(Q;P))$ is given by
$\bsd(\L(Q),\L(P))$. An element $(F_0 \prec F_1 \prec \cdots \prec
F_r , G)$ corresponds to the convex hull of the copies of $G \preceq
P$ within the $F_i$'s:
\begin{equation*}
  \conv \left( \widehat{F}_i + \varepsilon (G - \widehat{\kappa(F_i)})
    \right) = \conv \left( \widehat{F}_i - \varepsilon
      \widehat{\kappa(F_i)} \right) + \varepsilon G
\end{equation*}
This polytope is the image under an affine embedding of the product of
$G$ with an $r$-simplex.

The notion generalizes to the situation of two polyhedral complexes
with a realized order preserving, non-rank-increasing correspondence
between their face posets. We will use this in Section~\ref{sec:vf}
for the identity correspondence of $\L(T)$.
\begin{lemma} \label{lemma:mix}
  There is a coherent subdivision of $\dDS$ which is combinatorially
  isomorphic to the restriction to $|\Sigma|$ of the product
  subdivision $\bsd(S) \times T$.
\end{lemma}
\begin{proof}
  Here, $\Dv$ is a Minkowski summand of $\DSv$. Let $\psi$ be a
  piecewise linear concave function with domains of linearity given by
  $\bsd(\DSv;\Dv)$. Now $\nu$ induces another piecewise linear
  (non-concave) function $\tilde{\nu}$ on $\DSv$, which is $\nu$ in
  $G$-direction on $(F_0 \prec F_2 \prec \cdots \prec F_r , G)$, and
  constant in $F_i$-direction.
  \begin{figure}[htbp]
    \begin{center}
      \input{mix.pstex_t}
      \caption{The function $\tilde{\nu}$ on $\conv \left[ 
          \begin{smallmatrix}
            -2&\ 1&\ 1&\ 2&\ 2\\
            \ 0&\ 3&-3&\ 2&-2\\
            -1&-1&-1&-1&-1
          \end{smallmatrix} \right] \prec \DSv$.}
    \end{center}
  \end{figure}
  The function $\psi$ is strictly concave wherever $\tilde{\nu}$ is
  non-concave, so that for large $N$, the function $N \psi +
  \tilde{\nu}$ will be concave. Its domains of linearity are products
  of simplices that correspond to simplices of $\bsd(\DSv) \cong
  \bsd(S)$ times simplices of $T$.
\end{proof}
Now Lemma~\ref{lemma:sphere} follows if we start from $S$ and
$\bsd(T)$ in stead of $S$ and $T$.


%% file: delta.pstex_t
\begin{picture}(0,0)%
\includegraphics{delta.pstex}%
\end{picture}%
\setlength{\unitlength}{592sp}%
\begingroup\makeatletter\ifx\SetFigFont\undefined%
\gdef\SetFigFont#1#2#3#4#5{%
  \reset@font\fontsize{#1}{#2pt}%
  \fontfamily{#3}\fontseries{#4}\fontshape{#5}%
  \selectfont}%
\fi\endgroup%
\begin{picture}(10125,9361)(901,-8375)
\put(3376,-8236){\makebox(0,0)[lb]{\smash{\SetFigFont{5}{6.0}{\rmdefault}{\mddefault}{\updefault}$-1$}}}
\put(11026,-5461){\makebox(0,0)[lb]{\smash{\SetFigFont{5}{6.0}{\rmdefault}{\mddefault}{\updefault}$-3$}}}
\put(8701,-1486){\makebox(0,0)[lb]{\smash{\SetFigFont{5}{6.0}{\rmdefault}{\mddefault}{\updefault}$-1$}}}
\put(901,-4111){\makebox(0,0)[rb]{\smash{\SetFigFont{5}{6.0}{\rmdefault}{\mddefault}{\updefault}$-3$}}}
\put(5776,614){\makebox(0,0)[lb]{\smash{\SetFigFont{5}{6.0}{\rmdefault}{\mddefault}{\updefault}$-1$}}}
\end{picture}

%% file: deltaV.pstex_t
\begin{picture}(0,0)%
\includegraphics{deltaV.pstex}%
\end{picture}%
\setlength{\unitlength}{592sp}%
\begingroup\makeatletter\ifx\SetFigFont\undefined%
\gdef\SetFigFont#1#2#3#4#5{%
  \reset@font\fontsize{#1}{#2pt}%
  \fontfamily{#3}\fontseries{#4}\fontshape{#5}%
  \selectfont}%
\fi\endgroup%
\begin{picture}(7725,8236)(976,-9650)
\put(8701,-3061){\makebox(0,0)[lb]{\smash{\SetFigFont{5}{6.0}{\rmdefault}{\mddefault}{\updefault}$-1$}}}
\put(3526,-6061){\makebox(0,0)[rb]{\smash{\SetFigFont{5}{6.0}{\rmdefault}{\mddefault}{\updefault}$-1$}}}
\put(3151,-1786){\makebox(0,0)[lb]{\smash{\SetFigFont{5}{6.0}{\rmdefault}{\mddefault}{\updefault}$-1$}}}
\put(5926,-5461){\makebox(0,0)[rb]{\smash{\SetFigFont{5}{6.0}{\rmdefault}{\mddefault}{\updefault}$-3$}}}
\put(976,-5311){\makebox(0,0)[rb]{\smash{\SetFigFont{5}{6.0}{\rmdefault}{\mddefault}{\updefault}$-3$}}}
\put(4501,-9511){\makebox(0,0)[lb]{\smash{\SetFigFont{5}{6.0}{\rmdefault}{\mddefault}{\updefault}$-1$}}}
\end{picture}

%% file: pull.pstex_t
\begin{picture}(0,0)%
\includegraphics{pull.pstex}%
\end{picture}%
\setlength{\unitlength}{1184sp}%
\begingroup\makeatletter\ifx\SetFigFont\undefined%
\gdef\SetFigFont#1#2#3#4#5{%
  \reset@font\fontsize{#1}{#2pt}%
  \fontfamily{#3}\fontseries{#4}\fontshape{#5}%
  \selectfont}%
\fi\endgroup%
\begin{picture}(13888,5594)(1157,-5836)
\end{picture}

%% file: mix.pstex_t
\begin{picture}(0,0)%
\includegraphics{mix.pstex}%
\end{picture}%
\setlength{\unitlength}{1184sp}%
\begingroup\makeatletter\ifx\SetFigFont\undefined%
\gdef\SetFigFont#1#2#3#4#5{%
  \reset@font\fontsize{#1}{#2pt}%
  \fontfamily{#3}\fontseries{#4}\fontshape{#5}%
  \selectfont}%
\fi\endgroup%
\begin{picture}(7500,5896)(2851,-11855)
\put(5701,-9361){\makebox(0,0)[rb]{\smash{\SetFigFont{7}{8.4}{\rmdefault}{\mddefault}{\updefault}$-3$}}}
\put(5026,-7786){\makebox(0,0)[rb]{\smash{\SetFigFont{7}{8.4}{\rmdefault}{\mddefault}{\updefault}$-1$}}}
\put(4351,-8536){\makebox(0,0)[rb]{\smash{\SetFigFont{7}{8.4}{\rmdefault}{\mddefault}{\updefault}$-3$}}}
\put(2851,-10036){\makebox(0,0)[rb]{\smash{\SetFigFont{7}{8.4}{\rmdefault}{\mddefault}{\updefault}$-3$}}}
\put(3226,-10561){\makebox(0,0)[rb]{\smash{\SetFigFont{7}{8.4}{\rmdefault}{\mddefault}{\updefault}$-3$}}}
\put(3751,-11161){\makebox(0,0)[rb]{\smash{\SetFigFont{7}{8.4}{\rmdefault}{\mddefault}{\updefault}$-1$}}}
\put(4426,-11761){\makebox(0,0)[rb]{\smash{\SetFigFont{7}{8.4}{\rmdefault}{\mddefault}{\updefault}$-1$}}}
\put(6376,-11761){\makebox(0,0)[lb]{\smash{\SetFigFont{7}{8.4}{\rmdefault}{\mddefault}{\updefault}$-1$}}}
\put(6601,-10411){\makebox(0,0)[rb]{\smash{\SetFigFont{7}{8.4}{\rmdefault}{\mddefault}{\updefault}$-1$}}}
\put(6601,-8536){\makebox(0,0)[lb]{\smash{\SetFigFont{7}{8.4}{\rmdefault}{\mddefault}{\updefault}$-1$}}}
\put(7501,-9361){\makebox(0,0)[lb]{\smash{\SetFigFont{7}{8.4}{\rmdefault}{\mddefault}{\updefault}$-3$}}}
\put(8851,-8536){\makebox(0,0)[lb]{\smash{\SetFigFont{7}{8.4}{\rmdefault}{\mddefault}{\updefault}$-3$}}}
\put(8176,-7861){\makebox(0,0)[lb]{\smash{\SetFigFont{7}{8.4}{\rmdefault}{\mddefault}{\updefault}$-1$}}}
\put(6451,-6211){\makebox(0,0)[lb]{\smash{\SetFigFont{7}{8.4}{\rmdefault}{\mddefault}{\updefault}$-1$}}}
\put(10351,-10036){\makebox(0,0)[lb]{\smash{\SetFigFont{7}{8.4}{\rmdefault}{\mddefault}{\updefault}$-3$}}}
\put(9976,-10486){\makebox(0,0)[lb]{\smash{\SetFigFont{7}{8.4}{\rmdefault}{\mddefault}{\updefault}$-3$}}}
\put(9376,-11086){\makebox(0,0)[lb]{\smash{\SetFigFont{7}{8.4}{\rmdefault}{\mddefault}{\updefault}$-1$}}}
\put(8851,-11761){\makebox(0,0)[lb]{\smash{\SetFigFont{7}{8.4}{\rmdefault}{\mddefault}{\updefault}$-1$}}}
\end{picture}

%% file: monodromy.tex
\label{sec:monodromy}
An integral affine structure on an $n$-dimensional manifold $Y$ is
given by a torsion-free flat connection on the tangent bundle with
holonomy contained in $\SL(n,\Z)$. Alternatively, it can be
given by a coordinate covering $\{U_\alpha,
(y^i_\alpha)_{i = 1 \ldots n}\}$ of $Y$, together with transition
maps $f_{\alpha \beta} \in \SL(n,\Z) \ltimes
\R^n$ on the non-empty overlaps $U_\alpha \cap U_\beta$ such that the
usual cocycle condition is satisfied.
This endows the tangent bundle with the well defined lattice
structure, generated by the $\{\partial/\partial y^i_\alpha\}$.
\begin{defn}
  We define an integral affine structure on $\smooth$, using the covering 
  $\U \cup \V$, the nerve of which is the bipartite (two-colored) graph
  $\Gamma$: $U_v \cap V_w$ is non-empty if and only if 
  $\<v,w\> = 1$.
  For a point $q \in U_v$ we identify the tangent space $T_q(\smooth)$
  and the lattice $T^\Z_q$ in it with the following codimension 1
  subspace and sublattice of the pair $(\R^d,\Z^d)$:
  \begin{equation*}
    T_q = \R^d_v = \{n \in \R^d \suchthat \<v,n\> = 0 \} , \qquad
    T_q^\Z = \Z^d_v = \{n \in \Z^d \suchthat \<v,n\> = 0 \}.
  \end{equation*}
  For a point $q \in V_w$ we identify the tangent space $T_q(\smooth)$
  and the lattice in it with the $(d-1)$-dimensional
  quotients
  \begin{equation*}
    T_q = \R^d/w , \qquad T_q^\Z = \Z^d/w.
  \end{equation*}
  On the overlap $U_v \cap V_w$, we define the transition map $f_{vw}: \R^d_v 
  \rightarrow \R^d/w$ to be the restriction to the subspace $\R^d_{v}$ of the 
  natural projection $\R^d\rightarrow \R^d/w$.
\end{defn}
These transition maps respect the integral structure:
$f_{vw} \in \operatorname{Hom} (\Z^d_v, \Z^d/w)$, and the condition
$\<v,w\> = 1$ ensures that $f_{vw}$ is an isomorphism. The cocycle condition 
for the graph-type covering is trivial.

The rest of this subsection is devoted to description of the
monodromy. We call a loop in $\Gamma$ {\it primary}, if it consists of
4 edges: $(v_0, w_0),\ (w_0, v_1),\ (v_1, w_1),\ (w_1, v_0)$, for some
pair of edges $\{v_0, v_1\} \in S$, $\{w_0, w_1\} \in T$. We denote
such a loop by $(v_0w_0v_1w_1)$, and think of it as an element of the
fundamental group $\pi_1(\smooth) \cong \pi_1(\Gamma)$ with a base
point in $U_{v_0}$.
\begin{lemma}
  The monodromy transformation $T(v_0w_0v_1w_1) \colon \Z^d_{v_0}
  \rightarrow \Z^d_{v_0}$ along the loop $(v_0w_0v_1w_1)$ is given by
  $T(v_0w_0v_1w_1) (n) = n + \<v_1,n\>(w_1-w_0)$.
\end{lemma}
\begin{proof}
  Observe that $\<v_i,w_j\>=1,\ i,j=0,1$. Hence, if $n\in \Z^d_{v_0}$, 
  then $n + \<v_1,n\>(w_1-w_0) \in \Z^d_{v_0}$. 
  Now put $n' := n - \<v_1,n\> w_0 \in \Z^d_{v_1}$. Then
  $n' \equiv n \mod w_0$ together with $n' \equiv n +
  \<v_1,n\>(w_1-w_0) \mod w_1$ imply the desired formula.
\end{proof}
The fundamental group $\pi_1(\smooth) \cong \pi_1(\Gamma)$
is a free group whose abelianization $H_1(\Gamma)$ is generated (not
freely) by the primary loops. We apply the above calculation to
describe the local monodromy near a vertex in the discriminant locus.
\begin{defn}
  We will call a vertex $(\hat{\sigma},\hat{\tau})\in D$ to be of {\it type}
  $(k,l)$ if $\dim\sigma=k$ and $\dim\tau=l$. Let $\Star_\Sigma
  (\hat{\sigma},\hat{\tau})$ be its star neighborhood in $\Sigma$. Let
  $\vol(\sigma)$ denote the normalized volume of the simplex $\sigma$
  and same for $\vol(\tau)$.
\end{defn}
\begin{thm}
The fundamental group
$\pi_1(\Star_\Sigma (\hat{\sigma},\hat{\tau})\backslash D)$ is a free group
with $k\cdot l$ generators. In a suitable basis the monodromy
\begin{equation*}
  T \colon \pi_1(\Star_\Sigma (\hat{\sigma},\hat{\tau}) \backslash D)
  \rightarrow \SL(d-1,\Z)
\end{equation*}
can be represented (faithfully on the abelianization of $\pi_1$) by
an index $\vol(\sigma)\cdot \vol(\tau)$ subgroup of the (abelian) group of
matrices in the form  
\begin{equation*}
  \left( \begin{array}{cccccc} 1& 0& \ldots & *& \ldots & *\\ 
      0& 1& \ldots & \vdots& \ddots& \vdots \\
      0& 0& \ddots & *& \ldots &*\\
      \vdots& \vdots& \ddots & \ddots& \vdots& \vdots\\
      0& 0& \ldots & 0& 1& 0\\
      0& 0& \ldots& 0& 0& 1\\
    \end{array}\right)
\end{equation*}
(the identity matrix plus an $l\times k$ block in the right upper corner). 
\end{thm} 
\begin{proof}
Observe that $\Star_\Sigma (\hat{\sigma},\hat{\tau})\backslash
D$ is homotopy equivalent to the subgraph
$\Gamma_{(\hat{\sigma},\hat{\tau})}\subset\Gamma$, a complete bipartite graph
with $v$-nodes labeled by the vertices of $\sigma$, and $w$-nodes labeled by
the vertices of $\tau$. A choice of $v_0\in\sigma$ and $w_0\in\tau$ determines
the spanning tree of $\Gamma_{(\hat{\sigma},\hat{\tau})}$: 
\begin{center} 
  \begin{bipartite}{2cm}{4cm}{2cm}{3mm}{2mm}   \leftnode[0]{$v_0$}
  \leftnode[1]{\vdots} \leftnode[i]{$v_i$} \leftnode[k]{$v_k$}    
  \rightnode[10]{$w_0$} \rightnode[11]{\vdots} \rightnode[1j]{$w_j$}    
  \rightnode[13]{\vdots} \rightnode[1l]{$w_l$}     \match{0}{10} \match{0}{11}
  \match{0}{1j} \match{0}{13} \match{0}{1l}    \match{1}{10} \match{i}{10}
  \match{k}{10}   \brush{\dottedline{3}} \match{i}{1j} \end{bipartite} 
\end{center} 
Each edge $(v_iw_j),\ i=1,\dots,k,\ j=1,\dots,l$, added to the
spanning tree defines a primary loop which is a generator in the fundamental
group of $\Gamma_{(\hat{\sigma},\hat{\tau})}$. And there are no relations. 

Since $\<v_i,w_j-w_0\>=0$, all $i=0,\dots,k,\ j=1,\dots,l$, we can
choose a basis $\{e_r\}$ of $\Z^d_{v_0}$ such that
\begin{gather*}
w_j-w_0\in\mathrm{Span}_\Z\<e_1,\dots,e_l\>,\ \text{ all } j=1,\dots,l\\
\<v_i,e_r\>=0, \ \text{ all } i,r=1,\dots,k.
\end{gather*}
Then, in this basis the monodromy along the primary loop:
$(v_0w_0v_iw_j)$
\begin{equation*}
  T(v_0w_0v_iw_j)(e_r) = e_r + \<v_i,e_r\>(w_j-w_0)
\end{equation*}
will have the desired form. 

Finally, the factor of $\vol\sigma$ in the index of the subgroup
reflects the fact  that it may be impossible to complete the
collection $\{w_j-w_0\}$ to a full basis of $\Z^d_{v_0}$. Similarly,
$\vol\tau$ measures the failure of $\{v_i\}$ to form a part of a basis
for the dual lattice.
\end{proof}
\begin{remark}
If both triangulation $T$ and $S$ are unimodular, then $\vol(\sigma)\cdot
\vol(\tau)=1$, and one can
choose a basis $\{e_1,\dots,e_{d-1}\}$ of $\Z^d_{v_0}$ such that $e_j=w_j-w_0$
and $\<v_i,e_r\>=\delta_{d-i,r}$. The monodromy matrices then
will be
\begin{equation*}
  T(v_0w_0v_iw_j) = \mathbbm{1} + E(j,d-i),
\end{equation*}
where $E(j,d-i)$ is the elementary matrix with 1 in $(j,d-i)$-th place and
0 elsewhere. 
\end{remark}
\subsection{A glimpse of mirror symmetry} 
The invariance of $\Sigma$ and the discriminant locus $D$ with respect to the duality
between triangulated polytopes $\Delta$ and $\Delta^\vee$ is manifest. The
vertices of $D$ change the type from $(k,l)$ to $(l,k)$. The nodes of the
graph $\Gamma$ also interchange the $U$- and $V$-types. And if integral
bases associated with vertices of $\Gamma$ are chosen to be dual to the
original ones, then the transition matrices $f_{vw}$ will be replaced by
the transpose inverses $(f_{vw}^t)^{-1}$. Thus, on the same manifold
$Y=\Sigma\backslash D$ the two dual to each other integral affine
structures are realized.

Let us introduce notations for the following tori:
\begin{equation*}
\T:=\R^d/\Z^d, \qquad \T_v:=(\R^d_v)/(\Z^d_v), \qquad \T/w:=(\R^d/w)/(\Z^d/w). 
\end{equation*}
For $\<v,w\>=1$, the transition isomorphism
$f_{vw}\in\mathrm{Hom} (\Z^d_v,\Z^d/w)$ induces an isomorphism of the
tori, which we will denote by the same symbol
\begin{equation*}
  f_{vw} \colon \T_v \rightarrow \T/w.
\end{equation*}
The torus fibration over $Y=\smooth$ is constructed as follows. Using the
affine integral structure on $Y$ one can choose a covariantly constant
(with respect to the $\SL_\Z$-connection) lattice $T^\Z Y$ in the
tangent bundle $TY$ and form the relative quotient $W\rightarrow Y$
with the fibers $W_q=T_qY/T^\Z_qY$. Thus, the fibers are $W_q=\T_v$
when $q\in U_v$, and $W_q=\T/w$ when $q\in V_w$, with the canonical
identifications $f_{vw}: \T_v\rightarrow\T/w$ for $q\in U_v\cap V_w$.

Let $N(D)\subset\Sigma$ be a regular neighborhood of the discriminant locus. Let 
$W^\epsilon\to\base$ denote  the torus fibration associated to the original
integral affine  structure restricted to the complement of $N(D)$ in
$\Sigma$. In the next  section we will show that the torus fibration
$W^\epsilon$ on $\base$  embeds differentiably into $H_s$  for sufficiently
large $s$. The dual torus fibration by symmetry embeds into the mirror 
hypersurface. This is the topological part of the mirror symmetry
statement in the version of  Kontsevich and Soibelman \cite{KS}.


%% file: family.tex
\label{sec:family}
We describe the standard construction how to obtain the family of
Calabi-Yau hypersurfaces from our input
data~\cite{BatyrevMirror}. Recall that $\lambda \in \Z^{\D \cap (\Z^d)^*}$
and $\nu \in \Z^{\Dv \cap \Z^d}$ induce central triangulations $\{0\}
\join S$ of $\D$ and $\{0\} \join T$ of $\Dv$. That is, $\lambda$ and
$\nu$ lie in the interior of the secondary cone of the respective
triangulation~\cite{GKZ}.

We use $\lambda$ to define a (complex) one-parameter family $H_s^\aff$ of
affine hypersurfaces in the complex torus $(\C \setminus \{0\})^d$,
and we use $\nu$ in order to compactify the torus and the
hypersurfaces in a projective toric variety. The affine hypersurfaces
are given by 
\begin{equation*}
  H_s^\aff := \{x \in (\C \setminus \{0\})^d \suchthat \sum_{m \in \D
    \cap (\Z^d)^*} a_m s^{\lambda(m)} x^m = 0 \}
\end{equation*}
The triangulation $\{0\} \join T$ of $\Dv$ induced by the function
$\nu:\Dv\cap\Z\rightarrow\Z$ defines a simplicial
subdivision of the normal fan to $\D$, which in
turn defines a projective toric variety $X_\DT$ that contains $(\C
\setminus \{0\})^d$ as a dense open subset~\cite{Fulton,Oda,Danilov}. 
Let $H_s$ be the closure of $H_s^\aff$ in $X_{\DT}$.  

The function $\nu$ determines the class of an ample line
bundle on $X_{\DT}$, hence a K\"ahler class $[\nu]\in H^2(X_{\DT},\Z)$.
There are several, more or less canonical, ways to define a
$\T$-invariant K\"ahler form on $X_{\DT}$ in the class $[\nu]$. One of
the possible constructions of a toric variety is via symplectic
reduction. In this case $X_\DT$ inherits the natural symplectic
structure (which is, in fact, K\"ahler) from the standard K\"ahler
form $\frac{\sqrt{-1}}{4\pi}\sum dz\wedge d\bar{z}$ on $\C^{\vertT}$
(cf. \cite{Guillemin}).
Alternatively, we can use the pullback of the Fubini-Study form
$\omega_{\mathrm FS}$ from a projective embedding of $X_\DT$. Though
$\DT$ is not necessarily an integral polytope, some $k$-multiple of it
certainly is. We can use the complete linear system given by $k\DT$ to
define the embedding $i:X_\DT\hookrightarrow
\C\PP^{(k\DT)\cap(\Z^d)^*-1}.$ 
The K\"ahler form in the class $[\nu]\in H^2(X_\DT,\Z)$ is then given by 
$\frac1k i^* \omega_{\mathrm FS}$.  
   
In a sense any such ``canonical'' form $\omega_0$ is unsatisfactory
because the metric on $H_s$ defined by restriction of $\omega_0$ to
$H_s$ is too far from being Ricci-flat. In the second part of
this paper we will describe a family of forms $\omega_s$ on $X_\DT$,
such that the induced metrics on $H_s$ approximate the Calabi-Yau
metrics as $s\to\infty$ much better. (See the Outlook section for more
details). 

\begin{remark}
The toric variety $X_\DT$ is simplicial but not necessarily smooth, it
may have quotient singularities. Then we can
understand the K\"ahler forms in the orbifold sense (cf., e.g., \cite{AGM}).
\end{remark}

According to \cite[Ch.~10]{GKZ}, the hypersurfaces given by equations
in the form $\sum b_m x^m=0$ are all diffeomorphic to each other (in
the orbifold sense) as long as the vector $(\log|b_m|)$
($=\lambda\cdot \log|s|+\log|a|$ in our case) lies in some parallel
translation of the $(S\join\{0\})$-secondary cone. So that the
properties of the family related to the smooth structure do not depend
on along which ray we approach the large complex structure point
($s\to\infty$). Thus, any vector $\lambda$ in this secondary cone will
determine the diffeomorphic torus fibration. For the same reason we
can set the coefficients $a_m=1$ in the defining equation without loss
of generality. On the other side, any choice of the K\"ahler class, as
long as it is in the right K\"ahler cone, also gives rise to the same
combinatorics.

The goal of the rest of this section is to exhibit a torus fibration
$\Hsm\to\base$ on a ``smooth'' part of $H_s$, for large enough $s$, and
show that it is the same as our model fibration $W^\epsilon\to\base$.


%% file: moment1.tex
\label{sec:amoebas}
\newcommand{\nc}{\operatorname{NC}}
Let $\Log_s:(\C\setminus\{0\})^d\rightarrow\R^d$ be the logarithmic
map with the base $|s|$:
\begin{equation*}
  \Log_s(x) := \frac{\log(|x|)}{\log|s|} = \left\{
    \frac{\log|x_1|}{\log|s|}, \dots, \frac{\log|x_d|}{\log|s|}
  \right\}.
\end{equation*}
\begin{defn}(\cite[Ch.~6]{GKZ})
  The {\it amoeba} associated to the family of affine hypersurfaces
  $H_s^{\aff}$  is the image of the log map:
  \begin{equation*}
    \am_s := \Log_s(H_s^\aff) .
  \end{equation*}
\end{defn}
The geometry of amoebas of affine hypersurfaces is a well developed
subject that originated in the work of Gelfand, Kapranov and Zelevinsky
\cite{GKZ}. We are going to review several useful facts about the
amoebas, most of which are contained in (or can be easily deduced
from) a nice survey paper by Mikhalkin \cite{Mi1}.

The limiting behavior of amoebas as $s\to\infty$ can be described in
terms of the Legendre transform $L_\lambda:\R^d\rightarrow\R$ of the
vector $\lambda$:
\begin{equation*}
  L_\lambda(n) = \max_{m \in \D \cap (\Z^d)^*} \{\<m,n\> + \lambda(m)\}.
\end{equation*} 
\begin{remark}
  In the literature, the Legendre transform is sometimes defined with
  a ``minus'' rather than a ``plus'' sign. Those references work with
  convex (not concave) $\lambda$.
\end{remark}
$L_\lambda(n)$ is a piecewise linear convex function.  Define the {\it
non-Archimedean amoeba} $\am_\infty\subset\R^d$ to be the corner
locus of $L_\lambda(n)$ (the set of points where $L_\lambda(n)$ is not
smooth). $\am_\infty$ is a rational polyhedral complex of dimension
$d-1$ (cf. \cite{Mi2}), which gives a cell decomposition
of $\R^d$.

\begin{center} 
  \begin{figure}[htb]
    \epsfig{file=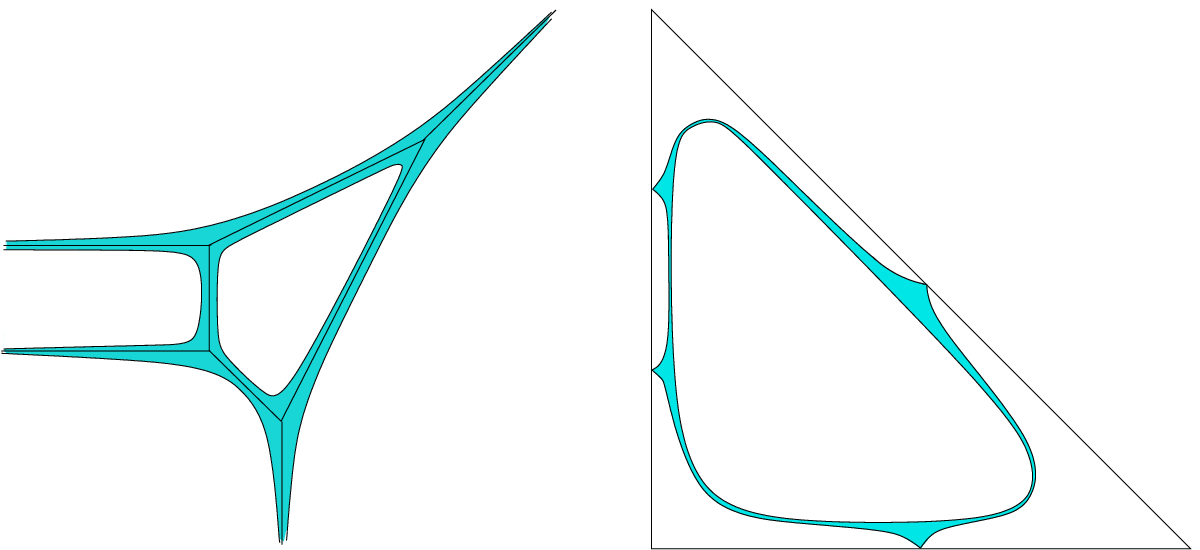} 
    \caption{\label{fig:am2}
      The
      affine amoeba $\am_s$ with the corresponding spine $\am_\infty$
      and its $\C\PP^2$-com\-pac\-ti\-fi\-ca\-tion $\overline{\am_s}$,
      for the family $H_s^{\aff} = \{s^2 + x^{-1}y^{-1} + sx^{-1} +
      sx^{-1}y^2 + sx^2y^{-1} = 0\}$.}
  \end{figure} 
\end{center}

\begin{lemma} \label{lemma:Q-cells}
  The decomposition of $\R^d$ by $\am_\infty$
  has the following description:
  \begin{enumerate}
  \item The cells are labeled by the simplices
    $\overline{\sigma} \in S \join \{0\}$.
  \item (The closure of) a cell $Q^\lambda_{\overline{\sigma}}$ is the
   Minkowski sum of the polytope and the cone:
    \begin{equation*}
      Q^\lambda_{\overline{\sigma}} = F^\vee_\sigma + \nc_\D(
      \overline{\sigma}).
    \end{equation*}
    Here $F^\vee_\sigma$ is the face of $\dDS$ dual to $\sigma =
    \overline{\sigma} \cap \dD \in S$, (where we set $F^\vee_\emptyset =
    \DSv$ for $\overline{\sigma} = \{0\}$), and
    $\nc_\D(F)$ is the normal cone to the face $F \prec \D$, (
    in particular $\nc_\D(F) = \{0\}$ for $F = \D$). 
    Thus, $Q^\lambda_{\overline{\sigma}}$ is unbounded if and only if
    $\overline{\sigma} = \sigma \in S$.
  \item In particular, the $d$-dimensional cells are labeled
    by the elements of
    $\vertS \cup \{0\}$. That is, there is a bounded central cell
    $Q^\lambda_{\{0\}} = \DSv$ and unbounded cells $Q^\lambda_{v}$, one for
    each vertex $v \in \vertS$ (see Fig.~\ref{fig:am2}).
  \end{enumerate}
\end{lemma}
\begin{proof}
  All statements follow easily from the definition of
  $L_\lambda$. Namely, the cells correspond to the subsets $I\subset
  \D \cap (\Z^d)^*$: the corresponding linear functions 
  $\<m,n\> +\lambda(m),\ m\in I$, saturate the maximum in $L_\lambda$.
  Since $\lambda$ is a concave function this can happen only
  if $I$ is a set of vertices of some simplex 
  $\overline{\sigma} \in S \join \{0\}$. This proves (1). 
  
  For (3) we notice that a $d$-cell is a domain of linearity of 
  $L_\lambda$, labeled by the vertex $m\in \vertS \cup \{0\}$ whose 
  corresponding linear 
  function $\<m,n\> +\lambda(m)$ is maximal. In particular,
  the central cell $Q^\lambda_{\{0\}}$ is the set of $n \in \R^d$, such that 
  the maximum is achieved by $\<\{0\},n\>+\lambda(0)$, i.e.
  \begin{equation*}
    \<m,n\> +\lambda(m)\leq\lambda(0), \ \text{ all } m\in\D\cap(\Z^d)^*,
  \end{equation*}
  which are exactly the defining inequalities for $\DSv$.
  
  More generally, a point $n$ is in (the closure of) 
  $Q^\lambda_{\overline{\sigma}}$ if and only if:
  \begin{gather*}
    \<m-v,n\> +\lambda(m)\leq \lambda(v)-\lambda(m), \ 
    \text{ all } v\in \mathrm{vert}\overline{\sigma},\ m\in\D\cap(\Z^d)^*,\\
    \text{ and }\ \<v_1,n\> +\lambda(v_1)=\<v_2,n\> +\lambda(v_2), \ 
    v_1,v_2\in \mathrm{vert}(\overline{\sigma}),
  \end{gather*} 
  which are exactly the defining inequalities for the polyhedron 
  $F^\vee_\sigma + \nc_\D (\overline{\sigma})$.
\end{proof}
The polyhedral complex $\am_\infty$ is also called the {\it spine} of
the amoeba $\am_s$ because of the following fact (cf. \cite{Mi2}):
\begin{prop} \label{prop:Hlimit}
  As $s\to\infty$ the amoebas $\am_s$ converge in the
  Hausdorff sense to the non-Archimedean amoeba $\am_\infty$.  
\end{prop}
\begin{proof}
[Idea of the proof]
  If we consider $s$ as a variable,
  we can think of the affine family $H_s^\aff$ as one hypersurface
  given by a single equation in $(\C\setminus\{0\})^{d+1}$. Then the rescaled
  amoeba $\log|s|\cdot\am_s$ sits inside the trace left by this extended
  $(d+1)$-dimensional amoeba in the horizontal hyperplane
  $\log|s|=const$. 
  And the result follows from \cite[Ch.~6, Prop.~1.9]{GKZ}.
\end{proof}
Given a $\T$-invariant K\"ahler form on $X_{\DT}$ in the class
$[\nu]$, we can consider the corresponding moment map
$\moment \colon X_{\DT} \to {\DT}$. In this case we can also define
the {\it compactified amoeba\/} $\overline{\am_s} := \moment(H_s)
\subset \DT$.

For the proof of Theorem~\ref{thm:main} we will need to introduce some
domains in $\R^d$, which are intimately connected with the amoebas and
the function $L_\lambda$. In some sense they are generalizations of
the cells induced by $\am_\infty$.

\begin{defn}
 For $I$ and $J$, two disjoint collections of integral points in
 $\D$, and  a real number $\epsilon\geq 0$, we define the (possibly empty) 
 polyhedron $\QeIJ$ in $\R^d$  by the conditions:
 \begin{gather*}
  \<m'',n\> + \lambda(m'') < \<m,n\> + \lambda(m) - \epsilon,\\
  \<m',n\> + \lambda(m') = \<m,n\> + \lambda(m),
 \end{gather*}
 for all $m,m'\in I, m''\not\in I\cup J$. We
 will abbreviate $\QeIJ$ by $Q^\lambda_I(\epsilon)$ when $J$ is empty. 
\end{defn}

\begin{defn}
  For $w \in \dDv \cap \Z^d$
  let $w^\perp$ be the set of integral points in $(\carrier_{\Dv}
  w)^\vee$, that is
  \begin{equation*}
    w^\perp = \{ m \in \D \cap (\Z^d)^* \suchthat \<m,w\> = 1 \}.
  \end{equation*}  
  Then, the {\it truncated} polytope $\D\setminus w^\perp$ is defined
  as the convex hull of integral points of $\D$ which are not in $w^\perp$.
\end{defn} 

Notice that for small $\epsilon$, $\QeIJ$ and $\QIJ$ are combinatorially equivalent.

\begin{lemma} \label{lemma:Qe-cells}
  In certain special cases of later interest we can describe $\QIJ$ as
  follows:
  \begin{enumerate}
  \item If $J=\emptyset$, then $Q^\lambda_I(0)$ is
    non-empty if and only if $I$ is the set of vertices of some simplex 
    $\overline{\sigma} \in S\join \{0\}$, in which case
    $Q^\lambda_I(0)=Q^\lambda_{\overline{\sigma}}$.
  \item $\Qv$ contains the relative interior of the facet
    $G_v\subset\DSv$.
  \item $\Qw$ is the Minkowski sum of the polytope $\DSv$
    and the normal cone to $\carrier_{\D\setminus w^\perp} \{ 0 \}$: 
    \begin{equation*}
      \Qw = \DSv + \nc_{\D \setminus w^\perp}( \{ 0 \}). 
    \end{equation*}
   In particular, it contains $\DSv+\cone(w)$.
  \end{enumerate}
\end{lemma}

\begin{proof}
 For (1) we notice that the set of defining inequalities of $Q^\lambda_I(0)$ is
 exactly the condition that the maximum in $L_\lambda$ is saturated by the
 linear functions $\<m,n\> +\lambda(m),\ m\in I$. 
 
 For (2), note that $Q^\lambda_{v\join\{0\}}$
 is defined by the same inequalities as $\Qv$, plus an extra 
 condition: $\<v,n\>+\lambda(v)=\lambda(0)$. 
 
 \begin{figure}[htbp]
    \begin{center}
      \includegraphics{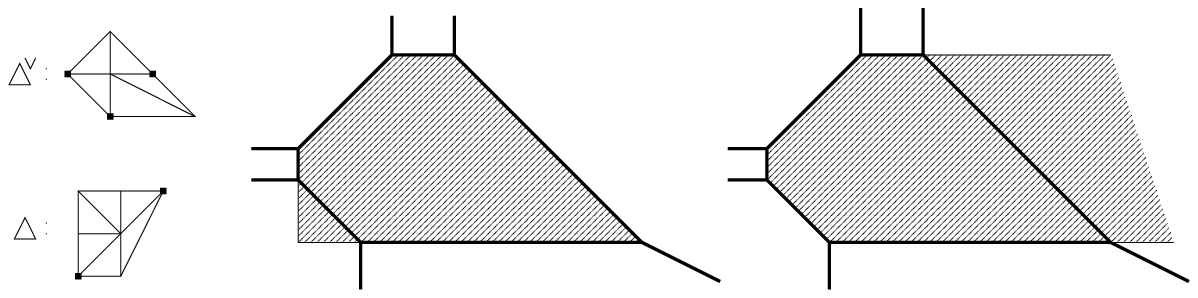}   
    \caption{\label{fig:qe2}
      Examples of $\Qv$ for 
      $v=\left[ 
        \begin{smallmatrix}
          -1\\-1
        \end{smallmatrix} \right],\
       \left[ 
        \begin{smallmatrix}
          1\\1
        \end{smallmatrix} \right]$
      in $\dD$.}
 \end{center}
 \end{figure} 
  
 For (3) we can study the Legendre transform of the restriction of $\lambda$
 to the truncated polytope $\D \setminus w^\perp$ (which is still a concave
 function). From this point of view, the Minkowski sum in (3) is in complete
 analogy with (2) of Lemma~\ref{lemma:Q-cells}. 
  \begin{figure}[htbp]
    \begin{center}
      \includegraphics{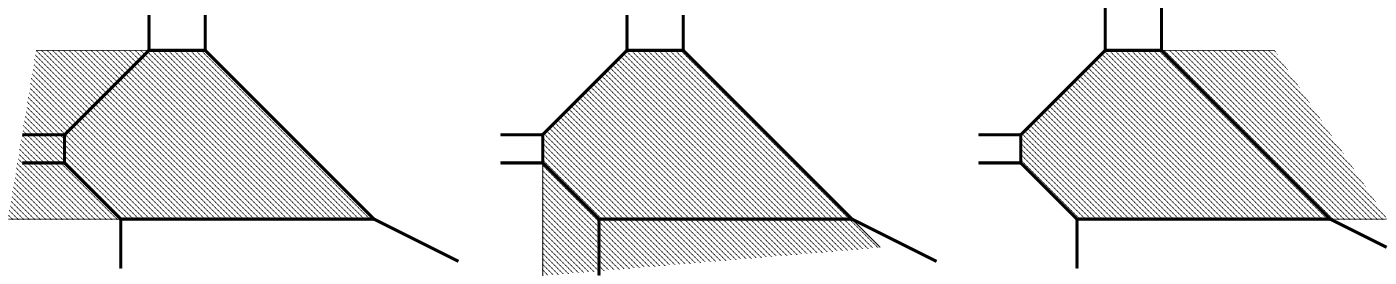}
    \caption{\label{fig:qe3}
      Examples of $\Qw$ for 
      $w=\left[ 
        \begin{smallmatrix}
          -1\\0
        \end{smallmatrix} \right],\  
      \left[ 
        \begin{smallmatrix}
          0\\-1
        \end{smallmatrix} \right],\
      \left[ 
        \begin{smallmatrix}
          1\\0
        \end{smallmatrix} \right]$      
      in $\dDS$.}
   \end{center}
 \end{figure} 
 
  Note that in the process of truncation we removed all integral 
 points of $\D$ with $\<m,w\>=1$. Hence, the remaining ones satisfy
 $\<m,w\>\leq 0$ which implies that $\{0\}$ is on the boundary of 
 $\D \setminus w^\perp$, and $w$ is in (the boundary of) the normal cone 
 $\nc_{\D \setminus w^\perp}(\{0\})$.
\end{proof}


%% file: moment2.tex
\label{sec:vf}
\newcommand{\dd}{\mathfrak{d}^\delta}
We will exhibit a vector field on $\R^d \setminus \ball$ with values
in $\dDS \subset \R^d$. Its integral curves yield the desired
foliation $\F$. Recall from Lemma~\ref{lemma:mix} that there is a
subdivision of $\dDS$ which is combinatorially isomorphic to $\bsd(S)
\times T$, restricted to $|\Sigma|$. In this context, the
projection $p_2$ can be thought of as defined $\dDS \rightarrow \dDv$,
and can be interpreted as a vector field on $\dDS$. We will deform
$p_2$, and extend the deformed vector field to $\R^d \setminus
\ball$.

For the deformation part, we use that the faces of $\bsd(T,T)$
are parameterized by pairs $(\tau, $ {\tiny $\mathfrak T$} $= \tau_0
\prec \ldots \prec \tau_r) \in T \times \bsd(T)$ with $\tau \preceq
\tau_0$ (cf.\ \S~\ref{sec:bsd}).
\begin{defn}
  In the $\delta$-realization of $\bsd(T,T)$, the face of $\bsd(T,T)$
  which corresponds to $(\tau, ${\tiny $\mathfrak T$}$)$ is the
  Minkowski sum $\delta \tau + (1-\delta)${\tiny $\mathfrak T$}.

  The $\delta$-realization of $\bsd(T,T)$ yields a cellular map
  $\dd \colon \bsd(T,T) \rightarrow T$, which maps $(\tau, $
  {\tiny $\mathfrak T$}$)$ to $\tau$. (I.e., the small copy of $\tau$
  in itself is stretched to full size, and the space between the small
  copies is collapsed into the smallest face.
\end{defn}
\begin{lemma} \label{lemma:const}
  If $\tau_1 \prec \tau_2$ are simplices of $T$, then $\dd$ is
  invariant with respect to $\widehat{\tau}_1 - \widehat{\tau}_2$ on
  the region $\bigcup_{\theta \in [0,1-\delta]} \theta
  \widehat{\tau_2} + (1-\theta) \tau_1$. In particular, $\dd$ is
  constant equal to $w$ on the whole star of $w \in \vertT$ in
  $\bsd(T,T)$.
\end{lemma}
\begin{figure}[htbp]
  \begin{center}
    \input{invariant.pstex_t}
    \caption{The map $\dd$ is ($\widehat{\tau}_1 -
      \widehat{\tau}_2$)-invariant in the shaded region.}
    \label{fig:delta}
  \end{center}
\end{figure}
\begin{lemma} \label{lemma:deform}
  Let $K$ be a subcomplex of $\bsd(T)$, and let $N$ be a
  neighborhood of $K$. Then there is a $\delta$-realization of
  $\bsd(T,T)$ such that for each face {\tiny $\mathfrak T$} of
  $K$, the faces of $\bsd(T,T)$ which correspond to some $(\tau,
  $ {\tiny $\mathfrak T$}$)$ are contained in $N$.
\end{lemma}
We will apply Lemma~\ref{lemma:deform} for $K = p_2(\dV)$, and $N =
p_2(N_2(\dV))$, where $N_2(\dV)$ is a neighborhood of $\dV$ in
$\Sigma$.
\begin{proof}
  Choose $\delta$ small enough to ensure that
  {\tiny $\mathfrak T$}
  $+ \ \delta (\tau_0 - \widehat{\tau}_0) \subset N$ for every simplex
  {\tiny $\mathfrak T$} $ = ( \tau_0 \prec \ldots \prec \tau_r ) \in
  K$.
  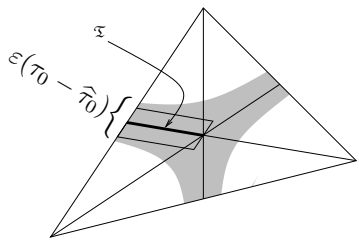
\begin{figure}[htbp]
    \begin{center}
      \input{deform.pstex_t}
      \caption{Choice of $\delta$ in the proof of
        Lemma~\ref{lemma:deform}.}
      \label{fig:deform}
    \end{center}
  \end{figure}

  Then the maximal cell which corresponds to $(\tau_0,$ {\tiny
    $\mathfrak T$}$)$ is given by
  \begin{equation*}
    \delta \tau_0 + (1-\delta) \text{\tiny $\mathfrak T$}
    \subseteq \text{\tiny $\mathfrak T$}
    + \delta (\tau_0 - \widehat{\tau}_0) \subset N .
  \end{equation*}
\end{proof}
Now we are ready to define the vector field $\X^\delta$ on $\dDS$ as
the composition $\dd p_2 \colon \dDS \rightarrow \dDv$.
In order to extend $\X^\delta$ to both sides of $\dDS$, we present a
polyhedral subdivision of a neighborhood of $\dDS$ whose trace on
$\dDS$ realizes the restriction of $\bsd(S) \times T$ to $|\Sigma|$.
\begin{remark}
  In \S~\ref{sec:bsd}, we were merely interested in the sphericity of
  $|\Sigma|$. We left open where to place the small copies of faces,
  and how small we wanted these copies to be. In the following we fix
  a realization of this subdivision which we will keep through the
  remainder of the article. In particular, $\epsilon$ is a fixed
  constant.
\end{remark}
Denote by $\lambda^\epsilon \in \R^{\D \cap (\Z^d)^*}$ the vector
given by $\lambda^\epsilon(0) = \lambda(0)$, and $\lambda^\epsilon(v)
= \lambda(v) + \epsilon$ for $v \in \vertS$. Suppose that $\epsilon >
0$ is small enough to ensure that $\lambda$ and $\lambda^\epsilon$
induce the same triangulation.  Then $\Dv_{\lambda^\epsilon} \subset
\DSv$ are combinatorially equivalent. For a simplex {\tiny
  $\mathfrak{S}$} $\in \bsd(S)$, denote {\tiny
  $\mathfrak{S}$}$_{\epsilon}$ the corresponding simplex of
$\bsd(\Dv_{\lambda^\epsilon})$. (I.e., {\tiny $\mathfrak{S}$} $\subset
(\R^d)^*$, while {\tiny $\mathfrak{S}$}$_{\epsilon} \subset \R^d$.)
Given a simplex $\tau \in T$ with $\<${\tiny $\mathfrak{S}$}$, \tau\>
= 1$, we can form the Minkowski sum {\tiny $\mathfrak{S}$}$_\epsilon +
\R_{\ge \epsilon/2} \tau$. These fit together to form a complex of
(unbounded) polyhedra which subdivides $\R^d$ outside
$\Dv_{\lambda^{\epsilon/2}}$. (It actually refines the subdivision
into $Q^{\lambda^{\epsilon/2}}_{\sigma}$'s provided by the
non-Archimedean amoeba for $\lambda^{\epsilon/2}$.)
\begin{figure}[htbp]
  \begin{center}
    \input{extend.pstex_t}
    \caption{The polyhedral subdivision of $\R^d$ outside
      $\Dv_{\lambda^\epsilon/2}$.}
    \label{fig:extend}
  \end{center}
\end{figure}
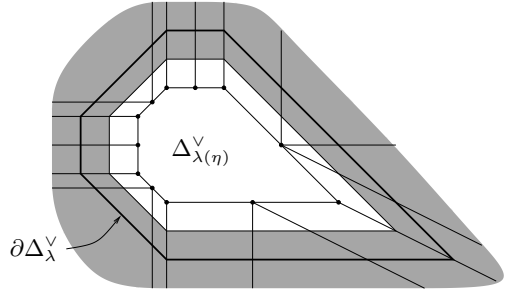
\begin{defn}
  For $0 < \delta < 1/2$, the vector field
  $\X^\delta \colon \R^d \setminus \Dv_{\lambda^{\epsilon/2}}
  \rightarrow \dDv \subset \R^d$ is the unique vector field which
  agrees with $\dd p_2$ on $\dDS$, and is invariant with respect to
  $\widehat{\tau}$ on the polyhedron {\tiny $\mathfrak{S}$}$_\epsilon
  + \R_{\ge \epsilon/2} \tau$.
\end{defn}
%
We need to argue that this determines a continuous vector field.
The problem may arise only when we try to assign a vector to a point $n
\in$ {\tiny $\mathfrak{S}$}$_\epsilon + \R_{\ge \epsilon/2} \tau$ such
that there are two points $n_1, n_2$ which already have a vector
assigned to them, so that both $n-n_1$ and $n-n_2$ are a multiple of
$\widehat{\tau}$.
Then $n_1 \in$ {\tiny $\mathfrak{S}$}$_\epsilon +
[\epsilon/2,\epsilon] \tau_1$ for some face $\tau_1 \prec \tau$, and
$\X^\delta(n_1)$ is defined as $\X^\delta(n_1')$, where $n_1' \in
\dDS$, and $n_1-n_1'$ is a multiple of $\widehat{\tau}_1$.
Furthermore,  $n_2 \in \dDS$, and $n_2 - n_1'$ is a multiple of
$\widehat{\tau} - \widehat{\tau}_1$. Here we use the assumption that
$\delta < 1/2$ to conclude that $\X^\delta(n_2) = \dd p_2(n_2) = \dd
p_2(n_1') = \X^\delta(n_1')$.
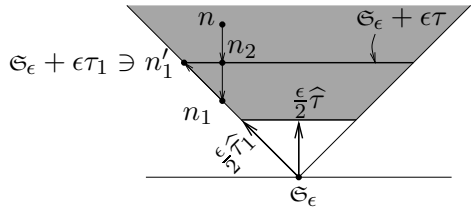
\begin{figure}[htbp]
  \begin{center}
    \input{welldef.pstex_t}
    \caption{$\X^\delta(n)$ is doubly defined: via $\X^\delta(n_1) =
      \X^\delta(n_1')$, and via $\X^\delta(n_2)$.}
    \label{fig:welldef}
  \end{center}
\end{figure}

The integral curves of $\X^\delta$ foliate $\R^d \setminus
\Dv_{\lambda^{\epsilon/2}}$. The following lemma summarizes the main
properties of $\X^\delta$ and the foliation $\F$. 
\begin{lemma} \label{lemma:amoeba}
  Given a neighborhood $N_2(\dV)$ of $\dV \subset \dDS \cong \Sigma$,
  there is a $\delta > 0$ such that
%
  \begin{enumerate}
  \item If $n \in Q^{\lambda^{\epsilon}}_v$, then
    $\<v,\X^\delta(n)\> = 1$.
  \item If $n \in V_w \setminus N_2(\dV)$, the flow line $\F_n$
    through $n$ is a straight line parallel to $w$ outside
    $\Dv_{\lambda^{\delta}}$.
  \end{enumerate}
\end{lemma}
\begin{figure}[htbp]
  \begin{center}
    \includegraphics[width=2in]{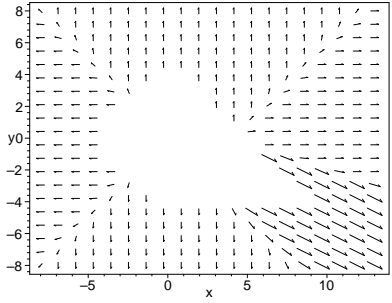}
    \caption{The vector field $\X^\delta$ for a $2$-dimensional example.}
    \label{fig:field2d}
  \end{center}
\end{figure}
\begin{proof}
  Choose $0 < \delta < 1/2$ so that the $2\delta$-realization of
  $\bsd(T,T)$ satisfies the conclusion of Lemma~\ref{lemma:deform} for
  $K = \dV$ and $N = p_2(N_2(\dV))$.

  By construction of $\X^\delta$, the set of values on one of
  the polyhedra {\tiny $\mathfrak{S}$}$_{\epsilon} + \R_{\ge
    \epsilon/2} \tau$ is contained in the set of values on its
  boundary which is contained in $\tau$. Statement (1) follows from
  $\<${\tiny $\mathfrak{S}$}$, \tau\> = 1$.

  For (2), let $n \in V_w \setminus N_2(\dV)$. Then $p_2(n)$ belongs
  to a cell $(${\tiny $\mathfrak{T}$}$, w)$ of the
  $2\delta$-realization of $\bsd(T,T)$, so that $\X^\delta(n) =
  w$. Also, say, $n \in \overline{U}_v$. If we parameterize $\F_n(t)$
  such that $\F_n(0) = n$ (and $\dot{\F}_n(t) = \X(\F_n(t))$), then,
  by (1), $\<v, \F_n(t)\> = \lambda(v) + t$. So $\F_n(t) - t
  \widehat{\tau}$ stays in the hyperplane $\<v,\cdot\>=1$, where $\tau
  = \carrier_T(${\tiny $\mathfrak{T}$}$)$.

  For $t > 0$, $\F_n(t) = n + t w$.
  For $-\delta < t < 0$, let $\ell \in (\R^d)^*$ be a linear
  functional which takes the values $0$ on $w$, and $1$ on the
  opposite side of {\tiny $\mathfrak{T}$}. Then $\ell(p_2(n)) <
  1-2\delta$, and $\frac{d}{dt} \ell(p_2(\F_n(t) - t \widehat{\tau}))
  = \ell(w-\widehat{\tau}) = 1$. So in this time range, $\F_n(t) = n +
  t w$ as well. For $t < -\delta$, $\F_n(t)$ belongs to
  $\Dv_{\lambda^{\delta}}$ by (1).
\end{proof}

\begin{remark}
  The foliation $\F$ can be, in fact, continued to the boundary of $\DT$
  (not smoothly at the $d-2$-skeleton of $\DT$) via the diffeomorphism
  $\moment \circ \Log_s^{-1}$ between $\R^d$ and the interior of
  $\DT$. So that it will induce a projection $X_\DT \setminus
  \Log_s^{-1}(\ball)\to\Sigma$. But to construct torus fibrations we
  will use only a part of this projection where it is clearly well
  defined. That is why we do not provide a proof for this more general
  statement here.
\end{remark}


%% file: invariant.pstex_t
\begin{picture}(0,0)%
\includegraphics{invariant.pstex}%
\end{picture}%
\setlength{\unitlength}{1579sp}%
\begingroup\makeatletter\ifx\SetFigFont\undefined%
\gdef\SetFigFont#1#2#3#4#5{%
  \reset@font\fontsize{#1}{#2pt}%
  \fontfamily{#3}\fontseries{#4}\fontshape{#5}%
  \selectfont}%
\fi\endgroup%
\begin{picture}(4845,3066)(2368,-3394)
\put(4951,-2311){\makebox(0,0)[lb]{\smash{\SetFigFont{10}{12.0}{\rmdefault}{\mddefault}{\updefault}{\color[rgb]{0,0,0}$\widehat{\tau}_2$}%
}}}
\put(3451,-1936){\makebox(0,0)[rb]{\smash{\SetFigFont{10}{12.0}{\rmdefault}{\mddefault}{\updefault}{\color[rgb]{0,0,0}$\widehat{\tau}_1$}%
}}}
\end{picture}

%% file: deform.pstex_t
\begin{picture}(0,0)%
\includegraphics{deform.pstex}%
\end{picture}%
\setlength{\unitlength}{1579sp}%
\begingroup\makeatletter\ifx\SetFigFont\undefined%
\gdef\SetFigFont#1#2#3#4#5{%
  \reset@font\fontsize{#1}{#2pt}%
  \fontfamily{#3}\fontseries{#4}\fontshape{#5}%
  \selectfont}%
\fi\endgroup%
\begin{picture}(5541,3838)(1672,-4187)
\put(3076,-811){\makebox(0,0)[lb]{\smash{\SetFigFont{6}{7.2}{\rmdefault}{\mddefault}{\updefault}{\color[rgb]{0,0,0}$\mathfrak{T}$}%
}}}
\put(1726,-1186){\rotatebox{330.0}{\makebox(0,0)[lb]{\smash{\SetFigFont{10}{12.0}{\rmdefault}{\mddefault}{\updefault}{\color[rgb]{0,0,0}$\varepsilon (\tau_0 - \widehat{\tau}_0) \Bigl\{$}%
}}}}
\end{picture}

%% file: extend.pstex_t
\begin{picture}(0,0)%
\includegraphics{extend.pstex}%
\end{picture}%
\setlength{\unitlength}{592sp}%
\begingroup\makeatletter\ifx\SetFigFont\undefined%
\gdef\SetFigFont#1#2#3#4#5{%
  \reset@font\fontsize{#1}{#2pt}%
  \fontfamily{#3}\fontseries{#4}\fontshape{#5}%
  \selectfont}%
\fi\endgroup%
\begin{picture}(18931,12024)(-11,-11173)
\put(7501,-5611){\makebox(0,0)[rb]{\smash{\SetFigFont{9}{10.8}{\rmdefault}{\mddefault}{\updefault}{\color[rgb]{0,0,0}$\Dv_{\lambda(\eta)}$}%
}}}
\put(301,-9811){\makebox(0,0)[rb]{\smash{\SetFigFont{9}{10.8}{\rmdefault}{\mddefault}{\updefault}{\color[rgb]{0,0,0}$\dDS$}%
}}}
\end{picture}

%% file: welldef.pstex_t
\begin{picture}(0,0)%
\includegraphics{welldef.pstex}%
\end{picture}%
\setlength{\unitlength}{789sp}%
\begingroup\makeatletter\ifx\SetFigFont\undefined%
\gdef\SetFigFont#1#2#3#4#5{%
  \reset@font\fontsize{#1}{#2pt}%
  \fontfamily{#3}\fontseries{#4}\fontshape{#5}%
  \selectfont}%
\fi\endgroup%
\begin{picture}(10824,6424)(589,-6149)
\put(3301,-3361){\makebox(0,0)[rb]{\smash{\SetFigFont{11}{13.2}{\rmdefault}{\mddefault}{\updefault}{\color[rgb]{0,0,0}$n_1$}%
}}}
\put(3376,-361){\makebox(0,0)[rb]{\smash{\SetFigFont{11}{13.2}{\rmdefault}{\mddefault}{\updefault}{\color[rgb]{0,0,0}$n$}%
}}}
\put(3751,-1261){\makebox(0,0)[lb]{\smash{\SetFigFont{11}{13.2}{\rmdefault}{\mddefault}{\updefault}{\color[rgb]{0,0,0}$n_2$}%
}}}
\put(5701,-5911){\makebox(0,0)[lb]{\smash{\SetFigFont{11}{13.2}{\rmdefault}{\mddefault}{\updefault}{\color[rgb]{0,0,0}{\tiny $\mathfrak{S}$}$_\epsilon$}%
}}}
\put(5776,-2986){\makebox(0,0)[lb]{\smash{\SetFigFont{11}{13.2}{\rmdefault}{\mddefault}{\updefault}{\color[rgb]{0,0,0}$\frac{\epsilon}{2} \widehat{\tau}$}%
}}}
\put(3676,-4936){\rotatebox{45.0}{\makebox(0,0)[lb]{\smash{\SetFigFont{11}{13.2}{\rmdefault}{\mddefault}{\updefault}{\color[rgb]{0,0,0}$\frac{\epsilon}{2} \widehat{\tau}_1$}%
}}}}
\put(7801,-361){\makebox(0,0)[lb]{\smash{\SetFigFont{11}{13.2}{\rmdefault}{\mddefault}{\updefault}{\color[rgb]{0,0,0}{\tiny $\mathfrak{S}$}$_\epsilon + \epsilon \tau$}%
}}}
\put(2101,-1786){\makebox(0,0)[rb]{\smash{\SetFigFont{11}{13.2}{\rmdefault}{\mddefault}{\updefault}{\color[rgb]{0,0,0}{\tiny $\mathfrak{S}$}$_\epsilon + \epsilon \tau_1 \ni n_1'$}%
}}}
\end{picture}

%% file: moment4.tex
\label{sec:fibration}
Using the foliation $\F$ we are going to define a
decomposition of the hypersurface $H_s=\Hsm\sqcup
H_s^{\mathrm{sing}}$, construct a torus fibration
$\Hsm\rightarrow\base$ and show that it is isomorphic to the fibration
$W^\epsilon\rightarrow\base$.

For any closed subset $J\subset\Sigma$ we will denote by
$X_s(J)\subset X_\DT$ the closure of $\Log_s^{-1}\left(\bigcup_{q\in
    J}\F_q\right)$ in $X_\DT$.
\begin{defn} 
  Let $N(D)$ be a regular neighborhood of $D$ in $\Sigma$. Then the
  {\it smooth} part of the hypersurface is $\Hsm := H_s \cap
  X_s(\base)$, and the rest $H_s^{\mathrm{sing}} := H_s \setminus
  \Hsm$ is {\it singular}.
\end{defn}
Since $D = \dU \cap \dV$, there exist regular neighborhoods $N_1(\dU)$
of $\dU$ and $N_2(\dV)$ of $\dV$ in $\Sigma$, such that $N(D) \supset
N_1(\dU) \cap N_2(\dV)$. This means that $\base$ can be covered by the
union of the closed sets:
\begin{equation*}
  \U^\epsilon=\{U^\epsilon_v\}=\{U_v\setminus 
  {N_1(\dU)}\}\ \text{ and }\ 
  \V^\delta=\{V^\delta_w\}=\{V_w\setminus {N_2(\dV)}\}.  
\end{equation*}
The amoebas $\am_s$, for a large enough $s$, all lie in
$\R^d\setminus\ball$. This means that $\F$ defines a projection
$\am_s\to\Sigma$ and, by composition with $\Log_s$, the projection
$H_s^\aff\to\Sigma$. Also $\am_s$ lie in $\R^d\setminus
Q^\lambda_v(\epsilon)$, for any $v\in\vertS$ and large $s$. Since the
unbounded ends of flow lines $\F_q$, for $q\in\Ue_v$, are in
$Q^\lambda_v(\epsilon)$ their closures do not contain any extra points
of the hypersurface:
\begin{equation*}
  H_s^\aff \cap \Log_s^{-1} \left( \bigcup_{q \in \Ue_v} \F_q \right)
  = H_s^\aff \cap X_s(\Ue_v) = H_s \cap X_s(\Ue_v).
\end{equation*}
Thus, the map  $H_s\cap X_s(\Ue_v)\to \Ue_v$ is well defined. On the
other hand, for two distinct points $q_1, q_2$ in $\Ve_w$ the
corresponding leaves are straight lines. Written in local coordinates
(see Lemma~\ref{lemma:w-fibers}) this implies that the sets
$X_s(\F_{q_1})$ and $X_s(\F_{q_2})$ are disjoint. Hence, the map
$H_s\cap X_s(\Ve_w)\to \Ve_w$ is well defined. Combined together we
have (for large enough $s$) the well defined projection
\begin{equation*}
  f_s \colon \Hsm \rightarrow \base, \qquad f_s(x) := q \
  \Leftrightarrow \ x \in X_s(q).
\end{equation*}
\begin{center}
 \begin{figure}[htb] 
  \epsfig{file=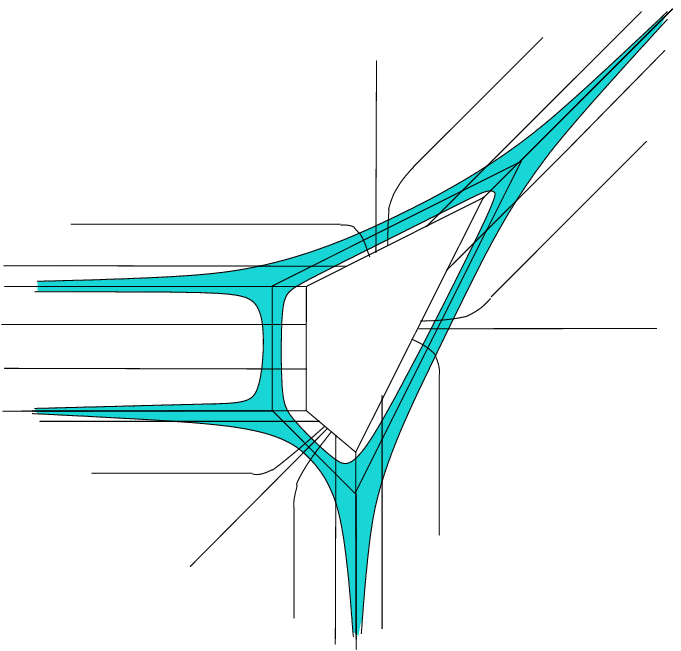} 
   \caption{The foliation $\F$ of $\R^d\setminus\ball$ induces a
     foliation of the amoeba $\am_s$.}
 \end{figure} 
\end{center}
\begin{thm} \label{thm:main} 
  There exists a real number $s_0$, such that for any $s$ with
  $|s|\geq s_0$, 
  \begin{equation*}
    f_s \colon \Hsm \rightarrow \base
  \end{equation*}
  is a torus fibration
  isomorphic to $\We\rightarrow\base$.
\end{thm}
Before proving the theorem we need to make a comment about
smoothness. $\Hsm$ is missing all singular points (if any) of the
toric variety $X_{\DT}$, which are all in the moment map preimage of
the $(d-2)$-skeleton of $\DT$ (see Lemma~\ref{lemma:w-fibers}). On the
other hand, $\base$ carries a canonical smooth structure induced by
the affine structure on $\smooth$. So given the topological fibration
$\Hsm \rightarrow \base$ of Theorem~\ref{thm:main}, standard
techniques apply to make it smooth. 

The strategy of proving Theorem~\ref{thm:main} will be as follows. First,
we analyze the map $f_s$ in $\Ue_v$ and $\Ve_w$ for every $v \in \vertS,
w \in \vertT$. Then the proof of the theorem can be completed in three
steps: we show that $f_s \colon \Hsm \to \base$ is a torus fibration
over the two kinds of covering patches, and then check that it has the
monodromy of our model.

Let $v \in \vertS$. For a fixed $s$ we consider the $(\D
\cap(\Z^d)^*-2)$-parameter family of hypersurfaces $H_s^v(a)$ in
$X_s(\Ue_v)$:
\begin{equation*}
  s^{\lambda(0)} + s^{\lambda(v)} x^v + \sum\limits_{m\neq \{0\},v}
  a_m s^{\lambda(m)} x^m = 0, \quad 0 \leq a_m \leq 1.
\end{equation*}
\begin{lemma} \label{lemma:v-fibers}
  There exists $s_0$ such that whenever $|s|\geq s_0$, all $H_s^v(a)$ are    
  smooth and transversal to $X_s(q)$ for every $q\in\Ue_v$. 
\end{lemma}
\begin{proof}
  According to Proposition~\ref{prop:Hlimit} we can choose $s$ big
  enough so that the $\Log_s$-image of every hypersurface $H_s^v(a)$
  lies in the $\epsilon$-neighborhood of $\Ue_v$. Recall from
  Lemma~\ref{lemma:Qe-cells} that a small neighborhood of $\Ue_v$ lies
  in the domain $\Qev$. Thus we can assume that all hypersurfaces
  $H_s^v(a)$ lie entirely in $\Log_s^{-1}(\Qev)$.
  
  Whenever $\Log_s(x)\in\Qev$, we have
  \begin{equation*}
    \<m,\frac{\log|x|}{\log|s|}\> + \lambda(m) \leq \lambda(0) -
    \epsilon, \text{ for all } m \neq v,\{0\} 
  \end{equation*}
  or, equivalently, 
  \begin{equation*}
    |x^m s^\lambda(m)| \leq |s|^{-\epsilon} |s|^{\lambda(0)}.
  \end{equation*}
  This means that the values of all monomials $x^ms^{\lambda(m)},\ m\neq
  v,\{0\}$, for $x\in \Log_s^{-1}(\Qev)$, are (uniformly) bounded by
  $|s|^{-\epsilon}|s|^{\lambda(0)}$. Note also, that their
  log-derivatives are bounded by $C|s|^{-\epsilon}|s|^{\lambda(0)}$,
  some constant $C\geq0$, since
  \begin{equation*}
    x \frac{\partial}{\partial x} (a_m s^{\lambda(m)}x^m ) = m \cdot a_m
    s^{\lambda(m)} x^m.
  \end{equation*}
  For any basis $\{e_i\}$ of $(\Z^d)^*$, the functions
  $y_i=x^{e_i}$ give affine coordinates on $(\Cstar)^d$. We choose $e_1=-v$,
   multiply the equations of the
  hypersurfaces in our family $H_s^v(a)$ by $y_1=x^{-v}$, and look for
  critical points:
  \begin{multline*}
    \frac{\partial}{\partial y_1}\bigl(y_1s^{\lambda(0)} +s^{\lambda(v)}
    +y_1\sum\limits_{m\neq \{0\},v}a_m s^{\lambda(m)}x^m \bigr) \\
    = s^{\lambda(0)}+\bigl(1+y_1\frac{\partial}{\partial y_1}\bigr)
    \sum\limits_{m\neq \{0\},v} a_m s^{\lambda(m)}x^m
    = s^{\lambda(0)} (1+O(|s|^{-\epsilon}))\neq 0,
  \end{multline*}
  for large enough $s$. Thus, there are no critical points, hence every
  member of our family $H_s^v(a)$ is smooth.
  
  Finally, note that $\bigcup_{q\in\Ue_v}\F_q$ is in $Q_v^{\lambda^\epsilon}$,
  but Lemma~\ref{lemma:amoeba} asserts that the vectors
   $\xi\in\X$ in $Q_v^{\lambda^\epsilon}$ satisfy 
  $\<v,\xi\>=1$. Thus, for any point of intersection $H_s^v(a)\cap X_s(q)$ the
  corresponding tangent vector to $X_s(q)$ has the form:
  $$ 
  \bar{\xi}= y_1\frac{\partial}{\partial y_1}+
  \alpha_2  y_2 \frac{\partial}{\partial y_2}+ 
  \dots+\alpha_d  y_d\frac{\partial}{\partial y_d}.
  $$ 
  Differentiating the defining equation for 
  $H_s^v(a)$ with respect to $\bar\xi$ gives:
  \begin{multline*}
    \bar\xi \bigl(y_1s^{\lambda(0)} +s^{\lambda(v)}
    +y_1\sum\limits_{m\neq \{0\},v}a_m s^{\lambda(m)}x^m \bigr)\\
    =y_1  s^{\lambda(0)} (1+O(|s|^{-\epsilon}))+ \sum\limits_{i=2}^{d} 
     \alpha_i y_i\frac{\partial}{\partial y_i} \bigl(y_1
    \sum\limits_{m\neq \{0\},v} a_m s^{\lambda(m)}x^m \bigr)\\
    =y_1  s^{\lambda(0)} (1+O(|s|^{-\epsilon}))+ \sum\limits_{i=2}^{d} 
    y_1  s^{\lambda(0)} O(|s|^{-\epsilon})
    =y_1 s^{\lambda(0)} (1+O(|s|^{-\epsilon}))\neq 0.
     \end{multline*} 
   Thus, we can conclude that $\bar\xi$ is transversal to the tangent planes to 
   $H_s^v(a)$, that is $X_s(q)$ is transversal to all $H_s^v(a)$. 
\end{proof}
\begin{remark}
  The estimates for the monomials in the lemma can be used to give
  another proof of the Hausdorff convergence in
  Proposition~\ref{prop:Hlimit}. Note that for any $\epsilon>0$ for
  large enough $s$, the monomial $x^v s^{\lambda(v)}$ become dominant in
  $\Log_s^{-1}(Q^\lambda_v(\epsilon))$. Hence the equation for $H_s^\aff$ 
  cannot have
  solutions in this domain. This means that the amoebas $\am_s$ are
  $\epsilon$-close to their spine $\am_\infty$.
\end{remark}
Now let $w\in\vertT$. Recall that $w^\perp$ is the set of integral points
in $(\carrier_{\Dv} w)^\vee$. For a fixed $s$ we consider the
  $(\D\cap(\Z^d)^*-w^\perp-1)$-parameter family of hypersurfaces:
  \begin{equation*}
    s^{\lambda(0)} + \sum\limits_{m\in G_w} s^{\lambda(m)} x^m +
    \sum\limits_{m \notin G_w \cup \{0\}} a_m s^{\lambda(m)} x^m = 0,
    \quad 0 \leq a_m \leq 1,
  \end{equation*}
and let $H_s^w(a)$ be its closure in $X_s(\Ve_w)$. Now we can repeat
the arguments of Lemma~\ref{lemma:v-fibers} to prove the analogous
statement for the family $H_s^w(a)$.
\begin{lemma} \label{lemma:w-fibers}
  There exists $s_0$ such that whenever $|s|\geq s_0$, all $H_s^w(a)$ are    
  smooth and transversal to $X_s(q)$ for every $q\in\Ve_w$.
\end{lemma}
\begin{proof}
  According to Proposition~\ref{prop:Hlimit} we can choose $s$ big enough so 
  that the $\Log_s$-image of the affine part of every hypersurface $H_s^w(a)$ 
  lies in
  the $\epsilon$-neighborhood of the Minkowski sum $\Ve_w+\cone(w)$. 
  Also, recall from Lemma~\ref{lemma:Qe-cells} that $\Ve_w+\cone(w)$ lies in
  the domain $\Qew$. Thus, we can assume that the affine parts of  all 
  hypersurfaces $H_s^w(a)$ lie in $\Log_s^{-1}(\Qew)$.
  
  We choose a basis $\{e_i\}$ of $(\Z^d)^*$ such that
  \begin{equation*}
    \<e_1,w\>=-1 \text{ and } \<e_i,w\>=0, i=2,\dots,d.
  \end{equation*}
  Then the affine coordinate functions $y_i=x^{e_i}$ can be extended
  (by allowing zero values for $y_1$) to the open part of the toric
  divisor $Z_w$ corresponding to the facet $F_w\subset\dDT$. Moreover, in these
  coordinates the preimage of each flow line $\F_q$ in $(\Cstar)^d$ is defined 
  by fixing the values of $|y_2|,\dots,|y_d|$, so that its closure $X_s(q)$ is 
  defined by the same equations, but allowing the zero value for $y_1$.  
   Hence, we can use $\{y_i\}$ as global coordinates on $X_s(\Ve_w)$.
 
  Multiplying the affine 
  equation of $H_s^w(a)$ by $y_1$ we note that the Laurent polynomial 
  $$ y_1\sum\limits_{m\in G_w}s^{\lambda(m)}x^m +
   y_1\sum\limits_{m\notin G_w\cup\{0\}}a_m s^{\lambda(m)}x^m
  $$
  has only positive powers of $y_1$, where as its first part
  $ P_1(y)=y_1\sum\limits_{m\in G_w}s^{\lambda(m)}x^m $
  is independent of $y_1$ at all.
  Thus, we get the global equation for the family in $X_s(\Ve_w)$.
  
  Now we can repeat the estimates for the monomials and their
  log-derivatives. Whenever $\Log_s(x)\in\Qew$, we have
  \begin{equation*}
    \<m,\frac{\log|x|}{\log|s|}\> + \lambda(m) \leq \lambda(0) -
    \epsilon, \text{ for all } m \notin G_w \cup \{0\}, 
  \end{equation*}
  or, equivalently,
  \begin{equation*}
    |x^ms^\lambda(m)|\leq |s|^{-\epsilon} |s|^{\lambda(0)}.
  \end{equation*}
  When written in the $y$-coordinates these estimates extends by continuity 
  from the affine part to the entire $X_s(\Ve_w)$. 
  
  To see that $H_s^w(a)$ has no critical points we differentiate its defining 
  equation  with respect to $y_1$:
  \begin{multline*}
    \frac{\partial}{\partial y_1} \bigl( y_1 s^{\lambda(0)} + y_1
    \sum\limits_{m \in G_w} s^{\lambda(m)} x^m + y_1 \sum\limits_{m
      \notin G_w \cup \{0\}} a_m s^{\lambda(m)} x^m \bigr) \\ 
    = s^{\lambda(0)}+\bigl(1+y_1\frac{\partial}{\partial y_1}\bigr) 
    \sum\limits_{m\notin G_w\cup\{0\}}a_m s^{\lambda(m)}x^m 
    = s^{\lambda(0)} (1+O(|s|^{-\epsilon}))\neq 0,
  \end{multline*}
  for large enough $s$.
 
  Finally, Lemma~\ref{lemma:amoeba} asserts that the vector field $\X$ in 
  $\bigcup_{q\in\Ve_w}\F_q$ is constant and equal 
  to $w$. It means that $\frac{\partial}{\partial y_1}$ is a tangent vector to
  $X_s(q)$, $q\in\Ve_w$, and it is transversal to $H_s^w(a)$ by the above 
  calculation.  
\end{proof}
\begin{proof}[Proof of Theorem~\ref{thm:main}]
  Note that if all $a_i=1$ in family $H_s^v(a)$, then we have the original
  equation of $H_s$. On the other hand, if all $a_i=0$, then  
  the family $H_s^v(a)$ degenerates to the hyperbola:
  \begin{equation*} 
    H_s^v(0) := \{ x \in X_s(\Ue_v) \suchthat s^{\lambda(0)} +
    s^{\lambda(v)} x^v = 0 \}.
  \end{equation*}
  Because $X_s(\F_q)$, $q\in\Ue_v$, intersect every $H_s^v(a)$
  transversally, the corresponding fibers $F_q:=H_s\cap X_s(\F_q)$ and
  $F_q^v:=H_s^v(0)\cap X_s(\F_q)$ are diffeomorphic.
 
  If $\theta=\{\theta_i\}$ denote the coordinates of the torus $\T$
  and $\theta_s$ is the phase of $s$, then the fiber $F_q^v$ of
  $H^v_s(0)$ is the torus
  \begin{equation*}
    F_q^v = \{ \theta \in \T \suchthat \<v,\theta\> +
    (\lambda(0)-\lambda(v)) \theta_s \equiv 0\ \mathrm{mod}\ 2\pi\},
  \end{equation*}
  which, for a fixed $s$, can be identified with the torus $\T_v$
  (though, see the remark below about monodromy as
  $\theta_s\mapsto\theta_s+2\pi$).

  Similarly, the fibers $F_q=H_s\cap X_s(\F_q)$ and $F_q^w:=H_s^w(0)\cap 
  X_s(\F_q)$ for $q\in\Ve_w$ are diffeomorphic. But $F_q^w$ can be naturally
  identified with the torus $\T/w$, which follows from writing the
  equation for $H_s^w(0)$ in the local coordinates $\{y_i\}$ from
  Lemma~\ref{lemma:w-fibers}:
  \begin{equation*}
    s^{\lambda(0)} y_1 + P_1(y_2,\dots,y_d) = 0,
  \end{equation*}
  where $P_1(y_2,\dots,y_d)$ is a Laurent polynomial independent of
  $y_1$. Restricting to the fiber $X_s(q)$ means fixing absolute values of
  $y_i,\ i=2,\dots,d$. A point on the torus $\T/w$ determines
  the phases of $y_i,\ i=2,\dots,d$. Once $y_i,\ i=2,\dots,d$,
  are fixed, there is a unique solution to the equation of $H_s^w(0)$.

  Thus, $f_s:\Hsm\rightarrow\base$ is a torus fibration. 
  The only thing left to check is that it has the correct monodromy. 

  Note that all diffeomorphisms $F_q\cong F^v_q$, $q\in\Ue_v$, and $F_q\cong 
  F^w_q$, $q\in\Ve_w$, 
  are deformation diffeomorphisms. Hence, the transitions maps between
  $\T_v$ and $\T/w$, for $q\in U_v \cap V_w$,
  are homotopic to the map $f_{vw}:\T_v\to\T/w$. But monodromy is a
  homotopy invariant, hence, it has to be equal to the one given by
  the maps $f_{vw}$. This completes the proof.
\end{proof}
\begin{remark}
  The same statement was proven in \cite{Zh} for regular hypersurfaces in 
  smooth 
  toric varieties using partition of unity arguments. This method can also be 
  applied in our situation since we do not touch the singular part of $X_\DT$ 
  at all. 
\end{remark}
\begin{remark}
  The fiber isomorphisms $F_p\cong\T_v$ depend on the value of the
  phase of $s$. If we go around a loop $s\mapsto se^{2\pi i}$, we
  won't come back to the original diffeomorphism $\Hsm\to\We$. Rather,
  it will be a composition with a generalized Dehn twist, namely, the
  diffeomorphism $\We\to\We$ which is the fiber wise shift by a section of 
  $\We\to\base$ (the tori are abelian
  groups). Such a section was explicitly written down in
  \cite{Zh}.
\end{remark}


%% file: outlook1.tex
\subsection{Combinatorics}
It would be interesting to know whether the subdivision of $|\Sigma|$
given by the $F \times F^\vee$ can be realized as the boundary of a
$d$-dimensional polytope (as the picture suggests). The dual of the
face lattice would be given by the lattice of intervals in the face
lattice of $\D$ (or of $\Dv$).

In fact, if there was a realization of the combinatorial
type of $\D$ with an identification $\R^d \cong (\R^d)^*$ such that
the normal cones of $F$ and $F^\vee$ intersect in their relative
interiors, then the Minkowski sum $\D + \Dv$ would do the trick. For
$d=3$, Koebe's Theorem~\cite{Thurston} guarantees such a realization.


%% file: outlook2.tex
\subsection{The Hausdorff convergence}
There is a natural family of complex structures $J_s$ on the model
torus fibration $h \colon W \to \smooth$. Namely, for a given $s$ we
can take $dy_i + \frac{\sqrt{-1}}{\log|s|} d \theta_i$ to be the
holomorphic 1-forms on $W$, where $\{y_i\}$ are the affine
coordinates on $\smooth$ and $\{\theta_i\}$ are the corresponding
coordinates on the torus fibers. Also, given a Riemannian metric $\sum
g_{ij} dy_i \otimes dy_j$ on $\smooth$, one can define the K\"ahler
metric on $W$ by $\omega_{ij}=\sum g_{ij}
(dy_i+\frac{\sqrt{-1}}{\log|s|}d\theta_i)\otimes
(dy_j-\frac{\sqrt{-1}}{\log|s|}d\theta_j)$. If, in addition, $g_{ij}$
satisfy the real Monge-Amp\`ere equation: $\det g_{ij} \equiv 1$, then
the induced metric on $W$ is Ricci-flat.

In the second part of the paper we will show that $\Hsm$ embeds
into $(W,J_s)$ ``almost'' holomorphically. Moreover, we will
construct a family of $\T$-invariant K\"ahler forms $\omega_s$ on
$X_\DT$ in the class $\frac{[\nu]}{\log|s|}$, such that the pairs
$(H_s,H_s^{\mathrm{sing}})$ with the induced metrics converge in the
Gromov-Hausdorff sense to the pair $(\Sigma,D)$. Here $\Sigma$ will
carry a compact metric space structure which will restrict to a
Riemannian metric on $\smooth$.

\subsection{Non-Archimedean geometry}
The family $H^\aff_s$ can be thought of as an affine algebraic
hypersurface in $(K^*)^d$ defined over a complete algebraically closed
field $K$ that contains $\C((s))$. Then, as the name suggests, the
image of $H_s^\aff$ under the valuation map $val:(K^*)^d\to R^d$ will
be the non-Archimedean amoeba $\am_\infty$ (cf. \cite{Kapranov}).

Interestingly, the potential functions for our K\"ahler forms $\omega_s$
will come from smoothing a convex piecewise linear function on
$\R^d$. This gives another evidence that there may be a reformulation of
mirror symmetry purely in non-Archimedean terms. There is a partial
understanding of this approach \cite{Rutgers2000}, but the entire
program still remains wide open.
